\documentclass[a4paper]{amsart}
\usepackage{amsfonts}
\usepackage{amssymb}
\usepackage{euscript}
\usepackage{xspace}
\usepackage{enumerate}

\newcommand{\blank}{{\llcorner\!\lrcorner}}
\newcommand{\Mp}[2][60]{${#2}$\hbox{-}\penalty#1\hspace*{0pt}}
\newcommand{\Mpn}[2][10000]{${#2}$\hbox{-}\penalty#1\hspace*{0pt}}

\newcommand{\Cstar}{C$^\ast$\hbox{-}\penalty10000\hspace*{0pt}}

\newcommand{\defeq}{\mathrel{:=}}

\newcommand{\CC}{\mathop{\mathrm{C}_\mathrm{c}}}
\newcommand{\CVI}{\mathop{\mathrm{C}_0}}
\newcommand{\BC}{\mathop{\mathrm{C}_\mathrm{b}}}

\newcommand{\NBC}{\mathop{\mathrm{C}}}

\newcommand{\C}{{\mathbb{C}}}
\newcommand{\R}{{\mathbb{R}}}

\newcommand{\Z}{{\mathbb{Z}}}
\newcommand{\N}{{\mathbb{N}}}
\newcommand{\Ztwo}{{\mathbb{Z}_2}}

\newcommand{\littleX}{\chi}
\newcommand{\midX}{\mathrm{X}}
\newcommand{\bigX}{\mathcal{X}}
\newcommand{\littleQ}{\mathrm{q}}
\newcommand{\bigQ}{\mathrm{Q}}
\newcommand{\Hilm}[1][E]{\EuScript{#1}}
\newcommand{\Hilsu}{\EuScript{H}}
\newcommand{\Hilsg}{\hat{\EuScript{H}}}

\newcommand{\ess}{{\mathrm{es}}}

\newcommand{\Adj}{\mathbb{L}}
\newcommand{\Comp}{\mathbb{K}}
\newcommand{\Mult}{\mathcal{M}}

\newcommand{\Mat}{\mathbb{M}}                   
\newcommand{\Cl}{\mathbb{C}\mathrm{l}}          
\newcommand{\Hils}[1][H]{\EuScript{#1}}         

\newcommand{\ket}[1]{\mathopen|{#1}\mathclose\rangle}
\newcommand{\bra}[1]{\mathopen\langle{#1}\mathclose|}
\newcommand{\5}[2]{\langle{{#1} \mid {#2}}\rangle}

\newcommand{\Ran}{\mathop\mathrm{Ran}}          

\newcommand{\Unse}[1]{#1^+}

\newcommand{\ID}{\mathrm{id}}

\newtheorem{theorem}{Theorem}[section]
\newtheorem{proposition}[theorem]{Proposition}
\newtheorem{corollary}[theorem]{Corollary}
\newtheorem{lemma}[theorem]{Lemma}

\theoremstyle{definition}

\theoremstyle{remark}

\newtheorem{remark}[theorem]{Remark}

\newcommand{\hot}{\mathbin{\hat{\otimes}}}

\newcommand{\subideal}{\mathrel\triangleleft}

\newcommand{\ev}{\mathrm{ev}}

\newcommand{\even}{\mathrm{even}}
\newcommand{\odd}{\mathrm{odd}}
\newcommand{\opp}{\mathrm{op}}

\newcommand{\congto}{\mathrel{\overset{\cong}{\longrightarrow}}}
\newcommand{\prto}{\twoheadrightarrow}
\newcommand{\injto}{\rightarrowtail}

\newcommand{\Cat}{\mathfrak{C}}
\newcommand{\Mor}{\operatorname{Mor}}
\newcommand{\GCalg}[1][G]{\text{$#1$-C$^\ast$}}
\newcommand{\SHo}[1][G]{[\text{$#1$-C$^\ast$}]_s}

\newcommand{\RKK}{\mathcal{R}KK}
\newcommand{\Ex}{\operatorname{Ex}}

\newcommand{\abs}[1]{\lvert#1\rvert}  
\newcommand{\norm}[1]{\lVert#1\rVert} 

\begin{document}

\author{Ralf Meyer}
\title[Equivariant KK-theory]{
  Equivariant Kasparov theory and generalized homomorphisms
}
\date{April 30, 2000}
\email{rameyer@math.uni-muenster.de}
\urladdr{http://wwwmath.uni-muenster.de/u/rameyer}
\address{ SFB 478---Geometrische Strukturen in der Mathematik\\
          Universit\"at M\"unster\\
          Hittorfstra\ss{}e~27\\
          48149 M\"unster\\
          Germany
  }
\begin{abstract}
  Let~$G$ be a locally compact group.  We describe elements of $KK^G (A,B)$ by
  equivariant homomorphisms, following Cuntz's treatment in the
  non-equivariant case.  This yields another proof for the universal property
  of $KK^G$: It is the universal split exact stable homotopy functor.
  
  To describe a Kasparov triple $(\Hilm, \phi, F)$ for $A,B$ by an equivariant
  homomorphism, we have to arrange for the Fredholm operator~$F$ to be
  equivariant.  This can be done if~$A$ is of the form $\Comp(L^2G) \otimes
  A'$ and more generally if the group action on~$A$ is proper in the sense of
  Exel and Rieffel.
\end{abstract}
\maketitle

\section{Introduction}
\label{sec:intro}

In this article, we carry over the description of Kasparov theory in terms of
generalized homomorphisms to the equivariant case.  Let us first recall the
well-known situation for Kasparov theory without group actions.

The existence and associativity of the Kasparov product mean that we can
define a category~$KK$ whose objects are the separable \Cstar{}algebras and
whose morphisms from $A$ to~$B$ are the elements of $KK (A,B)$.  In
\cite{cuntz:generalizedhom} and~\cite{cuntz:newlook}, Cuntz relates elements
of $KK(A,B)$ with trivially graded \Cstar{}algebras $A$ and~$B$ to ordinary
\Mpn{\ast}homomorphisms.  He defines a certain ideal $\littleQ A$ in the free
product $A \ast A$ and constructs a natural bijection between $KK (A,B)$ and
the set $[\littleQ A, \Comp \otimes B]$ of homotopy classes of
\Mpn{\ast}homomorphisms $\littleQ A \to \Comp \otimes B$, where~$\Comp$
denotes the algebra of compact operators on a separable Hilbert space.
Skandalis~\cite{skandalis:KKsurvey} remarks that we also have $KK (A,B) \cong
[\Comp \otimes \littleQ A, \Comp \otimes \littleQ B]$.  The Kasparov product
becomes simply the composition of \Mpn{\ast}homomorphisms in this picture.
Cuntz's description of $KK (A,B)$ is used by
Higson~\cite{higson:characterizeKK} to characterize Kasparov theory by a
universal property: The canonical functor from separable \Cstar{}algebras
to~$KK$ is the universal split exact stable homotopy functor.

For graded \Cstar{}algebras, Haag~\cite{haag:gradedKK} describes $KK (A,B)$ in
a similar way as the set of homotopy classes of grading preserving
\Mpn{\ast}homomorphisms from $\littleX A$ to $\hat{\Comp} \hot B$ for a
suitable graded \Cstar{}algebra $\littleX A$.  He shows that $KK (A,B) \cong
KK^{\Z_2} (\hat{S} \hot A,B)$, where $KK^{\Z_2}$ is the \Mpn{\Z_2}equivariant
Kasparov theory for trivially graded algebras and~$\hat{S}$ is $\CVI(\R)$
graded by reflection at the origin.  Furthermore, Haag identifies the Kasparov
product for graded \Cstar{}algebras in this setting~\cite{haag:algebraic}.

It is straightforward to carry over these results to $KK^G$ for a
\emph{compact} group~$G$.  However, new ideas are necessary if~$G$ is merely
locally compact.  The only result in that generality I am aware of is due to
Thomsen~\cite{thomsen:KKGuniversal}.  He shows that $KK^G$ can still be
characterized as the universal split exact stable homotopy functor.  However,
he does not obtain a description of $KK^G$ by equivariant
\Mpn{\ast}homomorphisms.

Let $A$ and~$B$ be \Mpn{G}\Cstar{}algebras.  Let~$\Comp$ be the algebra of
compact operators on the direct sum of infinitely many copies of $L^2G$.  We
would like to associate a \Mpn{G}equivariant \Mpn{\ast}homomorphism $\littleQ
A \to \Comp \otimes B$ (or $\littleX A \to \hat{\Comp} \hot B$ in the graded
case) to a Kasparov triple $(\Hilm, \phi, F)$ for $A,B$.  This may be
impossible for two reasons: The operator~$F$ need not be \Mpn{G}equivariant,
and there may be no \Mpn{G}equivariant embedding $\Hilm \subset L^2(G,
B)^\infty$.  It is surprisingly easy to overcome these problems: We just have
to replace~$A$ by $\Comp(L^2G) \otimes A$.  If the map $\phi \colon
\Comp(L^2G) \otimes A \to \Adj(\Hilm)$ is \emph{essential} in the sense that
$\phi(\Comp(L^2G) \otimes A) \cdot \Hilm$ is dense in~$\Hilm$, then $\Hilm =
L^2(G, \Hilm')$ for some Hilbert \Mp{B,G}module~$\Hilm'$.  Hence~$\Hilm$ can
be embedded in $L^2(G, B)^\infty$.  Moreover, the additional copy of~$L^2G$
gives us enough freedom to replace~$F$ by a \emph{compact perturbation}~$F'$
that is \Mpn{G}equivariant (Lemma~\ref{lem:equi_connection}).

Once we have that~$F$ is \Mpn{G}equivariant and $\Hilm \subset L^2(G,
B)^\infty$, we can proceed as in the non-equivariant case.  We show that we
get the same \Mpn{KK^G}groups if we restrict to Kasparov triples and
homotopies $(\Hilm, \phi, F)$ with a \Mpn{G}equivariant symmetry~$F$ and
$\Hilm \subset L^2(G,B)^\infty$ (Proposition~\ref{pro:special_aep}).
\emph{Symmetry} means that $F = F^\ast$ and $F^2 = 1$.  This yields a
bijection between $KK^G( \Comp(L^2G) \otimes A, B)$ and the set of homotopy
classes of \Mpn{G}equivariant \Mpn{\ast}homomorphisms $\littleQ (\Comp(L^2G)
\otimes A) \to \Comp \otimes B$ (Proposition~\ref{pro:KKGs_XQ}).  In addition,
we obtain an analogous statement for graded algebras and show that we may
tensor $\littleQ (\Comp(L^2G) \otimes A)$ with $\Comp(\Hils)$ for any
\Mpn{G}Hilbert space~$\Hils$.

The universal property of $KK^G$ for trivially graded separable
\Mpn{G}\Cstar{}algebras is an immediate consequence of this description of
$KK^G$ because $F(\littleQ A) \cong F(A)$ for any split exact stable homotopy
functor~$F$.  For graded algebras, we prove that
$$
KK^G(A, B) \cong KK^{G \times \Z_2} (\hat{S} \hot A, B),
$$
where $KK^G$ and $KK^{G \times \Z_2}$ denote the Kasparov theories for graded
\Mpn{G}\Cstar{}algebras and trivially graded \Mp{G \times
\Z_2}\Cstar{}algebras, respectively.  We describe the Kasparov product in this
setting.

In addition, we show that we can obtain $KK^G (A,B)$ using only Kasparov
triples $(\Hilm, \phi, F)$ with \Mpn{G}equivariant~$F$ and $\Hilm \subset
L^2(G, B)^\infty$ if~$A$ is \emph{proper} in the sense of
Exel~\cite{exel:MoritaSpectral} and Rieffel~\cite{rieffel:pre1}.  This notion
of properness is quite general and covers both algebras of the form
$\Comp(L^2G) \otimes A$ and the proper algebras
of~\cite{guentner-higson-trout}.

Our key result concerning proper group actions is that a countably generated
Hilbert \Mp{A,G}module~$\Hilm$ satisfies $\Hilm \oplus L^2(G,A)^\infty \cong
L^2(G,A)^\infty$ if and only if $\Comp(\Hilm)$ is a proper
\Mpn{G}\Cstar{}algebra.  This is not surprising in view of Rieffel's treatment
of square-integrable representations of groups on Hilbert
space~\cite{rieffel:pre1}.

\section{Notation and Conventions}
\label{sec:prep}

For the convenience of the reader, we recall the definitions of Hilbert
modules, Kasparov triples, and connections.  Moreover, we fix some notation.

Let~$G$ be a locally compact, \Mpn{\sigma}compact topological group.
Let~$dg$ be a left invariant Haar measure on~$G$ and let $L^2G = L^2(G,dg)$.
The \emph{left regular representation} of~$G$, defined by $\lambda_g(f) (g')
\defeq f(g^{-1} g')$ for $f \in L^2G$ and $g,g' \in G$, is a strongly
continuous unitary representation of~$G$ on~$L^2G$.  We always equip the
\Cstar{}algebra $\Comp(L^2G)$ of compact operators on~$L^2G$ with the
\Mpn{G}action induced by~$\lambda$.  That is, $\lambda_g(T) = \lambda_g \circ
T \circ \lambda_g^{-1}$ for all $g \in G$, $T \in \Comp(L^2G)$.

Let~$\Ztwo$ be the group with two elements and let $G_2 \defeq G \times
\Ztwo$.  A \emph{\Mpn{\Ztwo}graded \Mpn{G}\Cstar{}algebra} or, briefly,
\Mpn{G_2}\Cstar{}algebra is a \Cstar{}algebra with a strongly continuous
action of~$G_2$.  Recall that a grading is nothing but a \Mpn{\Ztwo}action.
We always write $\alpha_g$ and~$\beta_g$ for the actions of $g \in G_2$ on the
\Mpn{G_2}\Cstar{}algebras $A$ and~$B$, respectively.

Let $\Mat_2$ and~$\hat{\Mat}_2$ be the algebra of \Mp{2 \times 2}matrices with
the trivial grading and with the off-diagonal grading, respectively.  That is,
the off-diagonal terms in $\hat{\Mat}_2$ are odd.  Let $\Cl_1$ be the first
Clifford algebra, that is, the universal \Cstar{}algebra generated by an odd
symmetry.

\subsection{Hilbert modules}
\label{sec:Hilbert_modules}

Let~$B$ be a \Mpn{G_2}\Cstar{}algebras.  A \emph{\Mpn{\Ztwo}graded Hilbert
\Mp{B,G}module} or, briefly, \emph{Hilbert \Mp{B,G_2}module} is a Hilbert
\Mpn{B}module~$\Hilm$ with \Mpn{B}valued inner product $\5{\blank}{\blank}_B$
that is equipped with a strongly continuous linear action~$\gamma_g$ of~$G_2$
satisfying $ \gamma_g(\xi \cdot b) = \gamma_g(\xi) \beta_g(b) $ and $ \beta_g
(\5{\xi}{\eta}_B) = \5{\gamma_g\xi}{\gamma_g\eta}_B $ for all $g \in G_2$,
$\xi, \eta \in \Hilm$, $b \in B$.  We call~$\Hilm$ \emph{full} iff the linear
span of $\5{\Hilm}{\Hilm}_B$ is dense in~$B$.

We write $\Adj(\Hilm)$ and $\Comp(\Hilm)$ for the \Cstar{}algebras of
\emph{adjointable} and \emph{compact} operators on~$\Hilm$.  The latter is
generated by the \emph{rank one operators} $\ket{\xi}\bra{\eta}$ defined by
$$
\ket{\xi}\bra{\eta} (\zeta) \defeq \xi \cdot \5{\eta}{\zeta}_B
\qquad \text{for all $\xi,\eta,\zeta \in \Hilm$.}
$$
We always endow $\Adj(\Hilm)$ with the induced \Mp{G_2}action, $\gamma_g(T)
\defeq \gamma_g \circ T \circ \gamma_g^{-1}$.  This action is strongly
continuous on $\Comp (\Hilm)$ but usually not on $\Adj (\Hilm)$.  We call $T
\in \Adj (\Hilm)$ \emph{\Mpn{G}continuous} iff the map $g \mapsto
\gamma_g(T)$ is norm continuous.

We denote graded and ungraded spatial \emph{tensor products} of
\Mpn{\Ztwo}graded \Cstar{}algebras and Hilbert modules by $\hot$
and~$\otimes$, respectively.  If~$A$ is trivially graded, then there is no
difference between $A \hot B$ and $A \otimes B$.

\subsubsection{Standard Hilbert modules}
\label{sec:standard_modules}

Let $\ell^2(\N)$ be the separable Hilbert space with trivial \Mpn{G_2}action.
Let $\ell^2(\Ztwo\N)$ be the graded Hilbert space $\ell^2(\N)^\even \oplus
\ell^2(\N)^\odd$.  Let $L^2(G\N) \defeq L^2G \otimes \ell^2(\N)$ and
$L^2(G_2\N) \defeq L^2G \otimes \ell^2(\Ztwo\N)$.  We abbreviate $\Comp(G)
\defeq \Comp(L^2G)$, $\Comp(\N) \defeq \Comp\bigl(\ell^2(\N)\bigr)$, etc.
Moreover, we write $\Comp(\cdots) A$ instead of $\Comp(\cdots) \hot A \cong
\Comp(\cdots) \otimes A$.  Let~$\Hilm$ be a Hilbert \Mp{B,G_2}module.  Define
$\Hilm^\infty \defeq \ell^2(\N) \hot \Hilm$, $L^2(G, \Hilm) \defeq L^2G \hot
\Hilm$, and $L^2(G_2, \Hilm) \defeq L^2(G \times \Ztwo) \hot \Hilm$.  Let
$$
\Hilsg_B \defeq L^2(G_2, B)^\infty \cong L^2(G, B \oplus B^\opp)^\infty,
\qquad
\Hilsu_B \defeq L^2(G, B)^\infty.
$$
The Hilbert module~$\Hilsu_B$ is important only if~$B$ is trivially graded.

\subsubsection{Isometric embeddings}
\label{sec:iso_embed}

Let~$B$ be a \Mpn{G_2}\Cstar{}algebra and let $\Hilm$ and~$\Hilm[F]$ be
Hilbert \Mp{B,G_2}modules.  A map $\iota \colon \Hilm \to \Hilm[F]$ is called
an \emph{isometric embedding} iff it is a linear, \Mpn{G_2}equivariant
\Mpn{B}module map and satisfies $\5{\iota(\xi)} {\iota(\eta)}_B = \5{\xi}
{\eta}_B$ for all $\xi,\eta \in \Hilm$.  Hence~$\iota$ is injective and $\iota
(\Hilm) \subset \Hilm[F]$ is a closed \Mpn{G_2}invariant \Mpn{B}submodule.  We
do not require~$\iota$ to be adjointable.  This happens iff $\iota (\Hilm)$ is
\emph{complementable}, that is, $\iota (\Hilm) \oplus \iota (\Hilm)^\bot =
\Hilm[F]$.  We write $\Hilm \subset \Hilm[F]$ iff there is an isometric
embedding $\Hilm \to \Hilm[F]$.

\subsection{Hilbert bimodules}
\label{sec:Hilbert_bimodules}

Let $A$ and~$B$ be \Mpn{G_2}\Cstar{}algebras.  A \emph{\Mpn{\Ztwo}graded
  Hilbert \Mpn{A,B,G}bimodule} or, briefly, \emph{Hilbert
  \Mp{A,B,G_2}bimodule} is a Hilbert \Mp{B,G_2}module~$\Hilm$ with a
\Mpn{G_2}equivariant \Mpn{\ast}homomorphism $\phi \colon A \to \Adj(\Hilm)$.
We often use module notation for the action of~$A$ on~$\Hilm$, writing $a \xi$
instead of $\phi(a) (\xi)$.  The equivariance of~$\phi$ means that $\gamma_g
(a\xi) = \alpha_g (a) \gamma_g (\xi)$ for all $g \in G_2$, $a \in A$, $\xi \in
\Hilm$.

Let $A \cdot \Hilm \subset \Hilm$ be the subset of all elements of the form $a
\xi$ with $a \in A$, $\xi \in \Hilm$.  The Cohen-Hewitt factorization theorem
\cite{cohen:factorization}, \cite{hewitt-ross:harmonicII} implies that $A
\cdot \Hilm$ is a closed linear subspace.  We call~$\Hilm$ \emph{essential}
iff $A \cdot \Hilm = \Hilm$.  Let $\Mult(A)$ be the \emph{multiplier algebra}
of~$A$.  If~$\Hilm$ is essential, then there is a unique extension of~$\phi$
to a \Mpn{G_2}equivariant \Mpn{\ast}homomorphism $\phi \colon \Mult(A) \to
\Adj(\Hilm)$.  The extension is defined by $\phi(m) (a \cdot \xi) = (m\cdot a)
\cdot \xi$ for all $m \in \Mult(A)$, $a \in A$, $\xi \in \Hilm$.

If~$\Hilm_1$ is a Hilbert \Mp{A,B,G_2}module and~$\Hilm_2$ is a Hilbert
\Mp{B,C,G_2}module, then the tensor product $\Hilm_1 \hot_B \Hilm_2$ over~$B$
is defined as in~\cite{kasparov88:equivariantKK}.  It is a Hilbert
\Mp{A,C,G_2}bimodule.  If~$B$ acts on~$\Hilm_2$ via $\phi \colon B \to
\Adj(\Hilm_2)$, then we also use the more precise notation $\Hilm_1 \hot_\phi
\Hilm_2$ of~\cite{blackadar98:Ktheory}.

\subsection{Imprimitivity bimodules}
\label{sec:imprimitivity}

Let $A$ and~$B$ be \Mpn{G_2}\Cstar{}algebras.  A Hilbert \Mp{A,B,G_2}bimodule
$(\Hilm, \phi)$ is called an \emph{imprimitivity bimodule} iff it is full
and~$\phi$ is an isomorphism onto $\Comp(\Hilm)$ \cite{rieffel74:induced}.  We
call $A$ and~$B$ \emph{Morita-Rieffel equivalent} iff there is an
imprimitivity bimodule for them.  This is an equivalence relation.
Especially, if~$\Hilm$ is an imprimitivity bimodule for $A,B,G_2$, then there
is a \emph{dual} imprimitivity bimodule~$\Hilm^\ast$ for $B,A,G_2$.  It
satisfies $\Hilm^\ast \hot_A \Hilm \cong B$ as Hilbert \Mp{B,B,G_2}bimodules
and $\Hilm \hot_B \Hilm^\ast \cong A$ as Hilbert \Mp{A,A,G_2}bimodules.  A
concrete model for~$\Hilm^\ast$ is $\Comp(\Hilm,
B)$\label{dual_imprimitivity}.  The algebras $\Comp(\Hilm) \cong A$ and
$\Comp(B) \cong B$ operate on $\Comp(\Hilm, B)$ by composition.  The
\Mp{\Comp(\Hilm)}valued inner product is defined by $\5{T_1}{T_2} \defeq
T_1^\ast T_2$ for all $T_1, T_2 \in \Comp(\Hilm, B)$.

\subsection{Kasparov triples}
\label{sec:Kasparov_triples}

Let $A$ and~$B$ be \Mpn{\sigma}unital \Mpn{G_2}\Cstar{}algebras.  A
\emph{Kasparov triple for $A,B$} is a triple $(\Hilm,\phi,F)$, where $(\Hilm,
\phi)$ is a countably generated Hilbert \Mp{A,B,G_2}bimodule and $F \in
\Adj(\Hilm)$ is odd with respect to the grading and satisfies
\begin{equation}  \label{eq:KK_relations}
  [F,\phi(a)],\
  (1-F^2)\phi(a),\
  (F-F^\ast)\phi(a),\
  (\gamma_g(F) - F) \phi(a)
  \in \Comp(\Hilm)
\end{equation}
for all $a \in A$, $g \in G$.  The expression $[F,\phi(a)]$
in~\eqref{eq:KK_relations} is a \emph{graded commutator}.  In the following,
all commutators will be graded.  The Kasparov triple is called
\emph{degenerate} iff all the terms in~\eqref{eq:KK_relations} are zero.

Thomsen~\cite{thomsen:extensions} shows that~\eqref{eq:KK_relations} implies
that the operators $F \cdot \phi(a)$ are \Mpn{G}continuous for all $a \in A$.
Hence this additional requirement of Kasparov~\cite{kasparov88:equivariantKK}
is redundant.

Two Kasparov triples $(\Hilm_t, \phi_t, F_t)$, $t = 0,1$, are
\emph{unitarily equivalent} iff there is a \Mpn{G_2}equivariant unitary $U
\colon \Hilm_0 \to \Hilm_1$ with $\phi_1(a) U = U \phi_0(a)$ for all $a \in A$
and $F_1 U = U F_0$.  Up to unitary equivalence, Kasparov triples are
functorial for \Mpn{G_2}equivariant \Mpn{\ast}homomorphisms in both variables.
If $f \colon B_1 \to B_2$ is a \Mpn{G_2}equivariant \Mpn{\ast}homomorphism and
$(\Hilm, \phi, F)$ is a Kasparov triple for $A,B_1$, then
$$
f_\ast(\Hilm, \phi, F) \defeq
(\Hilm \hot_f B_2, \phi \hot 1, F \hot 1).
$$
Let $B[0,1] \defeq C([0,1]; B)$ with the pointwise action of~$G_2$ and let
$\ev_t \colon B[0,1] \to B$ be the evaluation homomorphism at $t \in [0,1]$.
A \emph{homotopy} between two Kasparov triples $T_0$ and~$T_1$ is a Kasparov
triple $\bar{T} = (\bar{\Hilm}, \bar{\phi}, \bar{F})$ for $A, B[0,1]$ such
that $\bar{T}|_t \defeq (\ev_t)_\ast (\bar{\Hilm}, \bar{\phi}, \bar{F})$ is
unitarily equivalent to~$T_t$ for $t = 0,1$.  The \emph{Kasparov group}
$KK^G(A,B)$ is defined as the set of homotopy classes of Kasparov triples for
$A,B$.

Let $(\Hilm, \phi, F)$ be a Kasparov triple for $A,B$.  We call $F' \in
\Adj(\Hilm)$ a \emph{compact perturbation} of~$F$ iff
$$
(F' - F)\phi(a) \in \Comp(\Hilm)
\quad \text{and} \quad
\phi(a)(F' - F) \in \Comp(\Hilm)
\qquad
\text{for all $a \in A$.}
$$
If~$F'$ is a compact perturbation of~$F$, then $(\Hilm,
\phi, F')$ is a Kasparov triple as well.  The triples $(\Hilm, \phi, F)$ and
$(\Hilm, \phi, F')$ are \emph{operator homotopic} via the obvious path $F_t
\defeq (1-t)F + tF'$, and therefore also homotopic.

\subsection{Connections}
\label{sec:connections}

Let~$\Hilm_1$ be a Hilbert \Mp{A,G_2}module and let~$\Hilm_2$ be a Hilbert
\Mp{A,B,G_2}bimodule.  Let $\Hilm_{12} \defeq \Hilm_1 \hot_A \Hilm_2$.  For
$\xi \in \Hilm_1$, define an adjointable operator $T_\xi \colon \Hilm_2 \to
\Hilm_{12}$ by $T_\xi(\eta) \defeq \xi \hot \eta$ and
$T_\xi^\ast( \eta \hot \zeta) \defeq \5{\xi}{\eta}_A \cdot \zeta$.  For $\xi
\in \Hilm_1$, $F_2 \in \Adj(\Hilm_2)$, and $F_{12} \in \Adj(\Hilm_{12})$,
define adjointable operators on $\Hilm_2 \oplus \Hilm_{12}$ by
$$
\tilde{T}_\xi \defeq
\begin{pmatrix} 0 & T_\xi^\ast \\ T_\xi & 0 
\end{pmatrix}
\qquad \text{and} \qquad
F_2 \oplus F_{12} \defeq
\begin{pmatrix} F_2 & 0 \\ 0 & F_{12} 
\end{pmatrix}.
$$
The operator~$F_{12}$ is called an \emph{\Mpn{F_2}connection} iff $[F_2
\oplus F_{12}, \tilde{T}_\xi] \in \Comp(\Hilm_2 \oplus \Hilm_{12})$ for all
$\xi \in \Hilm_1$.  Assume that $F_2$ and~$F_{12}$ are odd and self-adjoint
and denote the grading automorphism on~$\Hilm_1$ by~$\tau$.  Then~$F_{12}$ is
an \Mpn{F_2}connection iff
$$
F_{12} T_\xi - T_{\tau\xi} F_{2} \in \Comp(\Hilm_2, \Hilm_{12})
\qquad \text{for all $\xi \in \Hilm_1$}.
$$
We will freely use the standard properties of connections
\cite[18.3]{blackadar98:Ktheory}.

\section{Equivariant connections and special Kasparov triples}
\label{sec:equi_connections}

Let $A$ and~$B$ be \Mpn{\sigma}unital \Mpn{G_2}\Cstar{}algebras and
let~$\Hils$ be a separable \Mpn{G_2}Hilbert space.  A Kasparov triple $(\Hilm,
\phi, F)$ for $A,B$ is called \emph{\Mpn{\Hils}special} iff
\begin{enumerate}[(i)]
\item $F$ is a \Mpn{G}equivariant symmetry; and
\item $\Hils \hot \Hilm \subset \Hilsg_B$.
\end{enumerate}
An \emph{\Mpn{\Hils}special homotopy} is given by an \Mpn{\Hils}special
Kasparov triple for $A, B[0,1]$.  We let $KK^G_{s,\Hils} (A,B)$ be the set of
\Mpn{\Hils}special Kasparov triples modulo \Mpn{\Hils}special homotopy.  If
$\Hils = \C$, we omit the~$\Hils$ and talk about \emph{special triples},
\emph{special homotopies}, and $KK^G_s (A,B)$.  We are mostly interested in
the cases $\Hils = \C$ and $\Hils = L^2(G_2\N)$.  In the latter case, the
condition $\Hils \hot \Hilm \subset \Hilsg_B$ becomes tautological.  The
additional flexibility of choosing~$\Hils$ is useful in connection with
Proposition~\ref{pro:KKGs_XQ}.  A special triple is automatically
\Mpn{\Hils}special because $\Hils \hot \Hilsg_B \cong \Hilsg_B$.  Hence there
are canonical maps $KK^G_s (A,B) \to KK^G_{s, \Hils} \to KK^G(A,B)$.  Usually,
these maps fail to be isomorphisms.  For instance, if the unit element of
$KK^G(\C,\C)$ comes from an element of $KK^G_s (\C,\C)$, then~$G$ has to be
compact.

In this section, we show that $KK^G_s (A,B) \cong KK^G_{s,\Hils} (A,B) \cong
KK^G (A,B)$ if~$A$ has the property AE that is defined below.  We verify that
algebras of the form $\Comp(L^2G) A$ have this property.  In
Section~\ref{sec:proper}, we will see that proper algebras have property AE as
well.

\begin{lemma}  \label{lem:equi_connection}
  Let $A$ and~$B$ be \Mpn{\sigma}unital \Mpn{G_2}\Cstar{}algebras.  Let
  $(\Hilm, \phi, F)$ be an essential Kasparov triple for $A,B$.  Let $\Hilm'
  \defeq L^2(G,A) \hot_\phi \Hilm \cong L^2(G, \Hilm)$.

  There is a \Mpn{G}equivariant \Mpn{F}connection~$F'$ on~$\Hilm'$.  Even
  more, we can achieve that~$F'$ is a \Mpn{G}equivariant self-adjoint
  contraction.
\end{lemma}

\begin{proof}
  Let $\CC(G, \Hilm)$ be the space of continuous functions $G \to \Hilm$ with
  compact support.  The inner product $\5{f_1}{f_2}_B \defeq \int_G
  \5{f_1(g)}{f_2(g)}_B \,dg$ turns $\CC(G, \Hilm)$ into a pre-Hilbert
  \Mp{B}module.  Its completion is $L^2(G, \Hilm)$.  We have $L^2(G,A)
  \hot_\phi \Hilm \cong L^2(G, \Hilm)$ because~$\phi$ is essential.  We may
  assume that~$F$ is a self-adjoint contraction by
  \cite[17.4.3]{blackadar98:Ktheory}.  Define $F' \colon \CC(G, \Hilm) \to
  \CC(G, \Hilm)$ by
  $$
  (F'f)(g) = \gamma_g(F) f(g) =
  \gamma_g\bigl( F \gamma_g^{-1} f(g) \bigr) \qquad
  \text{for all $g \in G$, $f \in \CC(G, \Hilm)$.}
  $$
  It is straightforward to check that~$F'$ is \Mpn{G}equivariant and odd
  and extends to a self-adjoint contraction $F' \colon L^2(G, \Hilm) \to
  L^2(G, \Hilm)$.

  We claim that~$F'$ is an \Mpn{F}connection.  Denote the grading
  automorphisms on~$A$ and $L^2(G,A)$ by~$\tau$.  We have to check that $K
  \defeq T_\xi F - F' T_{\tau \xi} \in \Comp(\Hilm, \Hilm')$ for all $\xi \in
  L^2(G,A)$.  We may restrict to~$\xi$ of the form $\xi(g) = f(g) a$ with $f
  \in \CC(G)$, $a \in A$, because such elements span a dense subspace of
  $L^2(G,A)$.  We have
  $$
  (K \eta)(g) =
  f(g) \phi(a) F\eta - f(g) \gamma_g(F) \phi\tau(a) \eta =
  K_g(\eta)
  $$
  for all $\eta \in \Hilm$, where
  $$
  K_g \defeq f(g) \phi(a) F - f(g) \gamma_g(F) \phi\tau(a) =
  f(g) [\phi(a), F] + f(g) \bigl(F - \gamma_g(F) \bigr) \phi\tau(a).
  $$

  Since $(\Hilm, \phi, F)$ is a Kasparov triple and~$f$ has compact support,
  $K_g$ is a norm continuous compactly supported function $G \to
  \Comp(\Hilm)$.  Using a partition of unity, we can approximate the function
  $g \mapsto K_g$ uniformly by finite sums of functions $g \mapsto \psi(g) T$
  with $\psi \in \CC(G)$, $T \in \Comp(\Hilm)$.  Approximating~$T$ by sums of
  finite rank operators, we can approximate~$K$ in norm by finite sums of
  operators of the form $\eta \mapsto \psi \hot \ket{\xi} \bra{\zeta}
  (\eta)$.  Hence $K \in \Comp( \Hilm, \Hilm')$, so that~$F'$ is an
  \Mpn{F}connection.
\end{proof}

We say that a \Mpn{G_2}\Cstar{}algebra~$A$ has \emph{property AE} iff: For all
\Mpn{\sigma}unital \Mpn{G_2}\Cstar{}algebras~$B$ and all \emph{essential}
Kasparov triples $(\Hilm, \phi, F)$ for $A,B$, there is a \Mpn{G}equivariant
compact perturbation~$F'$ of~$F$ and there is an isometric embedding $\Hilm
\subset \Hilsg_B$.

The letters AE are an abbreviation for ``automatic equivariance''.

\begin{proposition}  \label{pro:stable_aep}
  Let $A$ and~$B$ be \Mpn{\sigma}unital \Mpn{G_2}\Cstar{}algebras and let
  $(\Hilm, \phi, F)$ be an \emph{essential} Kasparov triple for $\Comp(G) A,
  B$.  Then we can find a \Mpn{G}equivariant compact perturbation of~$F$ and
  an isomorphism
  $$
  \Hilm \oplus \Hilsg_B \cong \Hilsg_B
  $$
  of Hilbert \Mp{B,G_2}modules.

  Thus $\Comp(G) A \defeq \Comp(L^2G) \hot A$ has property AE.
\end{proposition}

\begin{proof}
  Let
  $$
  \psi \colon \Comp(G) A \congto \Comp\bigl( L^2(G,A) \bigr)
  $$
  be the canonical isomorphism.  Thus $(L^2(G,A), \psi)$ is an imprimitivity
  bimodule.  Let $(L^2(G,A)^\ast, \psi^\ast)$ be the corresponding dual
  imprimitivity bimodule.  That is, $L^2(G,A)^\ast$ is a Hilbert \Mp{\Comp(G)
  A, G}module and~$\psi^\ast$ is an isomorphism between $A$ and $\Comp(
  L^2(G,A)^\ast)$ such that
  $$
  L^2 (G, A) \hot_{\psi^\ast} L^2 (G, A)^\ast \cong
  \Comp(G) A
  $$
  as Hilbert \Mp{\Comp(G) A, \Comp(G) A,G_2}bimodules.  Let
  $$
  \Hilm_0 \defeq L^2(G, A)^\ast \hot_\phi \Hilm,
  \qquad
  \phi_0 \defeq \psi^\ast \hot 1 \colon A \to \Adj(\Hilm_0).
  $$
  Let $F_0 \in \Adj(\Hilm_0)$ be an \Mpn{F}connection.  Then $(\Hilm_0,
  \phi_0, F_0)$ is an essential Kasparov triple for $A,B$.  It is a Kasparov
  product of $(L^2(G,A)^\ast, \psi^\ast, 0)$ and $(\Hilm, \phi, F)$.
  
  Since~$\phi$ is essential, we have $\Comp(G) A \hot_\phi \Hilm
  \cong \Hilm$ and hence
  \begin{displaymath}
    \Hilm \cong
    L^2(G, A) \hot_{\psi^\ast} L^2(G, A)^\ast \hot_\phi \Hilm \cong
    L^2(G, A) \hot_{\phi_0} \Hilm_0
  \end{displaymath}
  as Hilbert \Mp{B,G_2}modules.  We have $\phi = \psi \hot 1 \colon \Comp(G)
  A \to \Adj(L^2(G,A) \hot_{\phi_0} \Hilm_0)$.

  By Lemma~\ref{lem:equi_connection}, there is a \Mpn{G}equivariant
  \Mpn{F_0}connection~$F'$ on~$\Hilm$.  The operator~$F'$ is an
  \Mpn{F}connection on $\Comp(G) A \hot_\phi \Hilm$
  by \cite[18.3.4.f]{blackadar98:Ktheory}.  Thus $F - F'$ is a
  \Mpn{0}connection.  This means that~$F'$ is a compact perturbation of~$F$ by
  \cite[18.3.2.c]{blackadar98:Ktheory}.  As a result, $F'$ is a
  \Mpn{G}equivariant compact perturbation of~$F$.

  The equivariant stabilization theorem \cite[Theorem
  2.5]{mingo-phillips:triviality} for the compact group~$\Ztwo$ yields
  $\Hilm_0 \oplus (B \oplus B^\opp)^\infty \cong (B \oplus B^\opp)^\infty$ as
  \Mpn{\Ztwo}graded Hilbert \Mpn{B}modules.  Hence
  $$
  \Hilm \oplus \Hilsg_B \cong
  L^2(G, \Hilm_0 \oplus (B \oplus B^\opp)^\infty) \cong
  L^2(G, (B \oplus B^\opp)^\infty) =
  \Hilsg_B
  $$
  as Hilbert \Mp{B,G_2}modules by \cite[Lemma
  2.3]{mingo-phillips:triviality}.  Thus $\Hilm \subset \Hilsg_B$.
\end{proof}

It is a well-known fact that any Kasparov triple is homotopic to an essential
triple \cite[18.3.6]{blackadar98:Ktheory}.  We need a more explicit
construction of the homotopy.

\begin{lemma}  \label{lem:make_essential}
  Let $A$ and~$B$ be \Mpn{\sigma}unital \Mpn{G_2}\Cstar{}algebras.  Let
  $(\Hilm, \phi, F)$ be a Kasparov triple for $A,B$.  Let $\Hilm_\ess \defeq
  \phi(A) \cdot \Hilm \cong A \hot_\phi \Hilm$ and define $\phi_\ess \colon A
  \to \Adj(\Hilm_\ess)$ by $\phi_\ess(a) = a \hot_\phi \ID_{{\Hilm}}$ for all
  $a \in A$.  Let~$F_\ess$ be an \Mpn{F}connection on~$\Hilm_\ess$.

  Then $(\Hilm_\ess, \phi_\ess, F_\ess)$ is a Kasparov triple.

  There is a canonical homotopy $(\bar{\Hilm}, \bar{\phi}, \bar{F})$ between
  $(\Hilm, \phi, F)$ and $(\Hilm_\ess, \phi_\ess, F_\ess)$.  We have
  $\bar{\Hilm} \subset (\Hilm \oplus \Hilm)[0,1]$.  The operator~$\bar{F}$ is
  a \Mpn{G}equivariant self-adjoint contraction if both $F$ and~$F_\ess$ are
  \Mpn{G}equivariant self-adjoint contractions.
\end{lemma}

\begin{proof}
  Define maps $\phi_{11} \colon A \to \Adj(\Hilm)$, $\phi_{12} \colon A \to
  \Adj(\Hilm_\ess, \Hilm)$, $\phi_{21} \colon A \to \Adj(\Hilm, \Hilm_\ess)$,
  $\phi_{22} \colon A \to \Adj(\Hilm_\ess)$ by $\phi_{ij}(a) \xi \defeq
  \phi(a) \xi$ for all~$\xi$ in the appropriate source $\Hilm$
  or~$\Hilm_\ess$.  These maps combine to a \Mpn{G_2}equivariant
  \Mpn{\ast}homomorphism
  $$
  \phi_\ast \defeq
  \begin{pmatrix}
    \phi_{11} & \phi_{12} \\ \phi_{21} & \phi_{22}
  \end{pmatrix} \colon
  \Mat_2(A) \to
  \begin{pmatrix}
    \Adj(\Hilm, \Hilm) &
    \Adj(\Hilm_\ess, \Hilm) \\
    \Adj(\Hilm, \Hilm_\ess) &
    \Adj(\Hilm_\ess, \Hilm_\ess)
  \end{pmatrix} =
  \Adj(\Hilm \oplus \Hilm_\ess).
  $$
  We claim that $T \defeq (\Hilm \oplus \Hilm_\ess, \phi_\ast, F \oplus
  F_\ess)$ is a Kasparov triple for $\Mat_2(A)$ and~$B$.  We have $\phi_{11} =
  \phi$ and $\phi_{22} = \phi_\ess$.  If $a \in A$, then $\phi_{21} (a)$ and
  $\phi_{12} (a^\ast)$ are the operators named $T_a$ and~$T_a^\ast$ in the
  definition of a connection in Section~\ref{sec:connections}.  Hence $[F
  \oplus F_\ess, \phi_\ast(x)] \in \Comp( \Hilm \oplus \Hilm_\ess)$ if~$x$ is
  off-diagonal.  Using $A \cdot A = A$, we can extend this to arbitrary $x \in
  \Mat_2(A)$.  The other conditions for a Kasparov triple like $(1- (F \oplus
  F_\ess)^2) \phi_\ast(x) \in \Comp(\Hilm \oplus \Hilm_\ess)$ follow easily
  from the standard properties of connections
  \cite[18.3.4]{blackadar98:Ktheory} if~$x$ is diagonal.  We can extend this
  to off-diagonal~$x$ using once again that $A \cdot A = A$.  Hence~$T$ is a
  Kasparov triple as asserted.

  Let $\iota_t \colon A \to \Mat_2(A)[0,1]$ be the rotation homotopy
  $$
  \iota_t(a) \defeq
  \begin{pmatrix}
    (1-t^2)a & t \sqrt{1-t^2} a \\
    t\sqrt{1-t^2} a & t^2 a
  \end{pmatrix}.
  $$
  We have
  $$
  (\iota_0)_\ast (T) = (\Hilm, \phi, F) \oplus (\Hilm_\ess, 0, F_\ess)
  \quad\text{and}\quad
  (\iota_1)_\ast (T) = (\Hilm, 0, F) \oplus (\Hilm_\ess, \phi_\ess, F_\ess).
  $$
  Thus up to degenerate triples $(\Hilm \oplus \Hilm_\ess, \phi_\ast \circ
  \iota_t, F \oplus F_\ess)$ is a homotopy between $(\Hilm, \phi, F)$ and
  $(\Hilm_\ess, \phi_\ess, F_\ess)$.  Using also the canonical homotopy
  between a degenerate triple and zero \cite[17.2.3]{blackadar98:Ktheory}, we
  obtain an explicit homotopy $(\bar{\Hilm}, \bar{\phi}, \bar{F})$ between
  $(\Hilm, \phi, F)$ and $(\Hilm_\ess, \phi_\ess, F_\ess)$.  Clearly,
  $\bar{\Hilm}$ and~$\bar{F}$ have the desired properties.
\end{proof}

\begin{proposition}  \label{pro:special_aep}
  Let $A$ and~$B$ be \Mpn{\sigma}unital \Mpn{G_2}\Cstar{}algebras and
  let~$\Hils$ be a separable \Mpn{G}Hilbert space.  Assume that~$A$ has
  property AE.  Then the canonical maps $KK^G_s (A,B) \to KK^G_{s,\Hils} (A,B)
  \to KK^G (A,B)$ are bijective.
  
  That is, any Kasparov triple for $A,B$ is homotopic to a special triple; if
  two \Mpn{\Hils}special triples are homotopic, then there is an
  \Mpn{\Hils}special homotopy between them.
\end{proposition}

\begin{proof}
  Let $(\Hilm, \phi, F)$ be a Kasparov triple for $A,B$.  We may replace the
  operator~$F_\ess$ in Lemma~\ref{lem:make_essential} by an arbitrary compact
  perturbation \cite[18.3.2.c]{blackadar98:Ktheory}.  Hence we may select a
  connection~$F_\ess$ that is a \Mpn{G}equivariant self-adjoint contraction by
  property~AE and \cite[17.4.2--3]{blackadar98:Ktheory}.  A standard trick
  \cite[17.6]{blackadar98:Ktheory} allows us replace $F_\ess$ by a symmetry.
  First add the degenerate triple $(\Hilm_\ess^{\opp}, 0, -F_\ess)$.  The
  operator
  $$
  \tilde{F} \defeq
  \begin{pmatrix}
    F_\ess & \sqrt{1-F_\ess^2} \\ \sqrt{1-F_\ess^2} & -F_\ess
  \end{pmatrix}
  \in \Adj(\Hilm_\ess \oplus \Hilm_\ess^\opp)
  $$
  is a \Mpn{G}equivariant symmetry and a compact perturbation of $F_\ess
  \oplus -F_\ess$.  Property AE implies that $\Hilm_\ess \oplus
  \Hilm_\ess^\opp \subset \Hilsg_B \oplus \Hilsg_B^\opp \cong \Hilsg_B$.  Thus
  $$
  \Psi(\Hilm, \phi, F) \defeq
  (\Hilm_\ess \oplus \Hilm_\ess^\opp, \phi_\ess \oplus 0, \tilde{F})
  $$
  is a special Kasparov triple that is homotopic to $(\Hilm_\ess, \phi_\ess,
  F_\ess)$ and hence to $(\Hilm, \phi, F)$ by Lemma~\ref{lem:make_essential}.
  The Kasparov triple $\Psi(\Hilm, \phi, F)$ is not quite well-defined
  because we have to choose a \Mpn{G}equivariant connection~$F_\ess$.
  Since~$F_\ess$ is determined uniquely up to a compact perturbation,
  $\Psi(\Hilm, \phi, F)$ is well-defined up to special homotopy.

  We have $\Psi \circ \Psi(\Hilm, \phi, F) = \Psi(\Hilm, \phi, F)$ because
  the essential part of $\Psi(\Hilm, \phi, F)$ is equal to $(\Hilm_\ess,
  \phi_\ess, F_\ess)$.  Assume that two special Kasparov triples of the form
  $\Psi(T_0)$ and $\Psi(T_1)$ are homotopic.  If we apply~$\Psi$ to a
  homotopy between them, we obtain a special homotopy between representatives
  of $\Psi \circ \Psi(T_j) = \Psi(T_j)$, $j = 0,1$.  Hence if Kasparov
  triples of the form $\Psi(T)$ are homotopic, then they are specially
  homotopic and a fortiori \Mpn{\Hils}specially homotopy.

  The proof will be finished if we show that if $T = (\Hilm, \phi, F)$ is an
  \Mpn{\Hils}special Kasparov triple, then there is an \Mpn{\Hils}special
  homotopy between $T$ and~$\Psi(T)$.

  By Lemma~\ref{lem:make_essential}, there is a homotopy $(\bar{\Hilm},
  \bar{\phi}, \bar{F})$ between $T$ and $(\Hilm_\ess, \phi_\ess, F_\ess)$ such
  that~$\bar{F}$ is a \Mpn{G}equivariant self-adjoint contraction and
  $\bar{\Hilm} \subset (\Hilm \oplus \Hilm)[0,1]$.  Thus $\Hils \hot
  \bar{\Hilm} \subset \Hilsg_B [0,1] = \Hilsg_{B[0,1]}$ because $\Hils \hot
  \Hilm \subset \Hilsg_B$.  Replacing~$\bar{F}$ by a symmetry as above, we
  obtain an \Mpn{\Hils}special homotopy between $T \oplus (\Hilm^\opp, 0, -F)$
  and $\Psi(T)$.  The canonical homotopy between~$T$ and $T \oplus
  (\Hilm^\opp, 0, -F)$ \cite[17.2.3]{blackadar98:Ktheory} is
  \Mpn{\Hils}special as well.
\end{proof}

\section{Isometric embeddings of Hilbert modules}
\label{sec:embed_Hilbert}

In this section, we provide some techniques to deal with not necessarily
adjointable embeddings of Hilbert modules.  Although the group action does not
create any additional difficulty here, we give complete proofs because the
corresponding arguments in \cite{cuntz:generalizedhom}, \cite{cuntz:newlook},
and~\cite{haag:gradedKK} are rather sketchy.

Let~$B$ be a \Mpn{G_2}\Cstar{}algebra and let $\Hilm$ and~$\Hilm[F]$ be
Hilbert \Mp{B,G_2}modules.  Let $\iota \colon \Hilm \to \Hilm[F]$ be an
isometric embedding as defined in Section~\ref{sec:iso_embed}.  Let
\begin{align*}
  \Adj_{\Hilm[F]}(\Hilm) & \defeq
  \{ T \in \Adj(\Hilm[F]) \mid
  \text{$T(\Hilm[F]) \subset \iota(\Hilm)$ and $T^\ast(\Hilm[F]) \subset
    \iota(\Hilm)$} \},
  \\
  \Comp_{\Hilm[F]}(\Hilm) & \defeq
  \Adj_{\Hilm[F]}(\Hilm) \cap \Comp(\Hilm[F]).
\end{align*}
Clearly, $\Adj_{\Hilm[F]} (\Hilm)$ and $\Comp_{\Hilm[F]} (\Hilm)$ are
\Cstar{}subalgebras of $\Adj(\Hilm[F])$.  If $T_1,T_2 \in \Adj_{\Hilm[F]}
(\Hilm)$, then $T_1 \Adj (\Hilm[F]) T_2 \subset \Adj_{\Hilm[F]} (\Hilm)$.
Thus $\Adj_{\Hilm[F]} (\Hilm)$ and $\Comp_{\Hilm[F]} (\Hilm)$ are hereditary
subalgebras of $\Adj (\Hilm[F])$.

\begin{lemma}  \label{lem:embedding_functorial}
  For $T \in \Adj_{\Hilm[F]} (\Hilm)$, define $\rho(T) \colon \Hilm \to \Hilm$
  by $\rho(T)(\xi) \defeq \iota^{-1} (T\iota\xi)$ for all $\xi \in \Hilm$.
  This yields a \Mpn{G_2}equivariant isometric \Mpn{\ast}homomorphism $\rho
  \colon \Adj_{\Hilm[F]} (\Hilm) \to \Adj (\Hilm)$.  Its restriction to
  $\Comp_{\Hilm[F]} (\Hilm)$ is an isomorphism onto $\Comp (\Hilm)$.

  Let $\Comp(\iota) \colon \Comp(\Hilm) \to \Comp_{\Hilm[F]} (\Hilm) \subset
  \Comp(\Hilm[F])$ be the inverse of $\rho|_{\Comp_{\Hilm[F]}(\Hilm)}$.  Then
  $\Comp(\iota)$ is the unique \Mpn{\ast}homomorphism satisfying
  \begin{equation}  \label{eq:comp_iota}
    \Comp(\iota) (\ket{\xi} \bra{\eta}) =
    \ket{\iota\xi} \bra{\iota\eta}
    \qquad \text{for all $\xi,\eta \in \Hilm$.}
  \end{equation}
\end{lemma}

\begin{proof}
  Clearly, $\rho(T)$ is adjointable for all $T \in \Adj_{\Hilm[F]}(\Hilm)$,
  with adjoint $\rho(T^\ast)$.  Thus~$\rho$ is a \Mpn{\ast}homomorphism
  $\Adj_{\Hilm[F]} (\Hilm) \to \Adj (\Hilm)$.  If $\rho(T) = 0$, then $T$
  vanishes on $\iota(\Hilm) \supset \Ran T^\ast$, so that $T \circ T^\ast =
  0$ and hence $T = 0$.  Thus~$\rho$ is isometric.  Since~$\rho$ is natural,
  it is \Mpn{G_2}equivariant.  If $\xi,\eta \in \Hilm$, then $\ket{\iota\xi}
  \bra{\iota\eta} \in \Comp_{\Hilm[F]} (\Hilm)$ and $\rho( \ket{\iota\xi}
  \bra{\iota\eta}) = \ket{\xi} \bra{\eta}$.  Thus $\rho\bigl(
  \Comp_{\Hilm[F]} (\Hilm) \bigr)$ contains $\Comp(\Hilm)$ and $\Comp(\iota)$
  satisfies~\eqref{eq:comp_iota}.
  
  It remains to show $\rho\bigl( \Comp_{\Hilm[F]} (\Hilm) \bigr) \subset
  \Comp (\Hilm)$.  It suffices to verify $\rho(TT^\ast) \in \Comp(\Hilm)$ for
  all $T \in \Comp_{\Hilm[F]} (\Hilm)$.  Evidently, $\rho(T \ket{\xi}
  \bra{\eta} T^\ast) = \rho( \ket{T\xi} \bra{T\eta})$ is a rank one operator
  for all $\xi, \eta \in \Hilm[F]$ because $T\xi, T\eta \in \iota (\Hilm)$.
  Therefore, $\rho(T u T^\ast) \in \Comp(\Hilm)$ for all $u \in
  \Comp(\Hilm[F])$.  If we let~$u$ run through an approximate unit for
  $\Comp(\Hilm[F])$, we get $\rho(TT^\ast) \in \Comp(\Hilm)$.
\end{proof}

By the way, if $\rho(T) = 1$, then $T^\ast T \colon \Hilm[F] \to
\iota(\Hilm)$ is a projection onto~$\iota(\Hilm)$, so that $\iota(\Hilm)$ is
complementable.  Hence~$\rho$ is surjective iff $\iota(\Hilm)$ is
complementable.  As an immediate consequence, we obtain the following result
of Combes and Zettl~\cite{combes-zettl}.

\begin{corollary}  \label{cor:hereditary_comp}
  Let~$B$ be a \Cstar{}algebra and\/~$\Hilm[F]$ a Hilbert \Mp{B}module.
  Let $H \subset \Comp (\Hilm[F])$ be a hereditary subalgebra.  Then $H =
  \Comp_{\Hilm[F]} (H \cdot \Hilm[F]) \cong \Comp (H \cdot \Hilm[F])$.

  Thus the hereditary subalgebras of\/ $\Comp (\Hilm[F])$ correspond
  bijectively to the not necessarily complementable Hilbert submodules
  of\/~$\Hilm[F]$.
\end{corollary}

\begin{proof}
  Since~$H$ is hereditary, $\ket{\xi} \bra{\eta} \in H$ for all $\xi, \eta
  \in H \cdot \Hilm[F]$.  By Lemma~\ref{lem:embedding_functorial}, these
  operators generate $\Comp_{\Hilm[F]} (H \cdot \Hilm[F])$.  Thus
  $\Comp_{\Hilm[F]} (H \cdot \Hilm[F]) \subset H$.  Obviously, $H \subset
  \Comp_{\Hilm[F]} (H \cdot \Hilm[F])$.
\end{proof}

Two isometric embeddings $\iota_0, \iota_1 \colon \Hilm \to \Hilm[F]$ are
\emph{homotopic} iff they can be connected by a continuous path of isometric
embeddings $\iota_t \colon \Hilm \to \Hilm[F]$, $t \in [0,1]$.  Such a
path~$\iota_t$ gives rise to an isometric embedding $h \colon \Hilm{}[0,1] \to
\Hilm[F][0,1]$, $(hf)(t) = \iota_t\bigl( f(t) \bigr)$.  The embedding~$h$
induces a map $\Comp(h) \colon \Comp(\Hilm)[0,1] \to \Comp(\Hilm[F])[0,1]$ by
Lemma~\ref{lem:embedding_functorial}.  Composing it with the inclusion
$\Comp(\Hilm) \to \Comp(\Hilm)[0,1]$ by constant functions, we obtain a
\Mpn{G_2}equivariant homotopy between $\Comp(\iota_0)$ and $\Comp(\iota_1)$.
As a result, homotopic isometric embeddings $\Hilm \to \Hilm[F]$ induce
homotopic \Mpn{\ast}homomorphisms $\Comp(\Hilm) \to \Comp(\Hilm[F])$.

\begin{lemma}  \label{lem:stabilization_homotopy}
  Let $B$ be a \Mpn{G_2}\Cstar{}algebra and let $\Hilm$ and\/~$\Hilm[F]$ be
  Hilbert \Mp{B,G_2}modules.  Then any two isometric embeddings\/ $\Hilm \to
  \Hilm[F]^\infty$ are homotopic.
\end{lemma}

\begin{proof}
  Let $\iota_0, \iota_1 \colon \Hilm \to \Hilm[F]^\infty$ be two isometric
  embeddings.  It is well-known that $\Hilm[F]^\infty \oplus \Hilm[F]^\infty
  \cong \Hilm[F]^\infty$ as Hilbert \Mp{B,G_2}modules, and that the
  inclusions of the direct summands $j_0,j_1 \colon \Hilm[F]^\infty \to
  \Hilm[F]^\infty$ are homotopic to the identity map.  These homotopies may
  be chosen \Mpn{G_2}equivariant.  Hence~$\iota_0$ is homotopic to $\iota'_0
  \defeq j_0 \circ \iota_0$ and~$\iota_1$ is homotopic to $\iota'_1 \defeq
  j_1 \circ \iota_1$.  By construction, $\iota'_0$ and~$\iota'_1$ have
  orthogonal ranges, that is, $\5{\iota'_0(\xi)} {\iota'_1(\eta)}_B = 0$ for
  all $\xi, \eta \in \Hilm$.  Hence $\iota'_t \defeq \sqrt{1-t^2} \iota'_0 +
  t\iota'_1 \colon \Hilm \to \Hilm[F]^\infty$ is an isometric embedding for
  all $t \in [0,1]$.  Thus $\iota'_0$ and~$\iota'_1$ are homotopic.
\end{proof}

The following lemma generalizes the observation of
Skandalis~\cite{skandalis:KKsurvey} that a degenerate Kasparov triple is
homotopic to zero.  It is also related to
\cite[Lemma~5.1]{cuntz:generalizedhom}.

\begin{lemma}  \label{lem:degenerate_homotopy}
  Let\/ $(\Hilm, \phi, F)$ be a Kasparov triple for $A,B$.  Let~$E$ be the
  \Cstar{}subalgebra of\/ $\Adj(\Hilm)$ generated by $\phi(A)$ and the
  operators $\gamma_g(F)$ for $g \in G$.  Let $J \subideal E$ be the smallest
  \Mpn{G}invariant ideal containing the operators
  $$
  [F,\phi(a)],\ 
  (1-F^2) \phi(a),\ 
  (F-F^\ast) \phi(a),\ 
  (\gamma_g(F) - F) \phi(a)
  $$
  for all $a \in A$, $g \in G$.  These are precisely the operators
  in~\eqref{eq:KK_relations} whose compactness (or vanishing) is required for
  a (degenerate) Kasparov triple.  Let $\Hilm' \defeq J \cdot \Hilm$.

  Then $\Hilm' \subset \Hilm$ is a closed, \Mpn{G_2}invariant submodule and
  $E (\Hilm') \subset \Hilm'$.  Hence restriction to~$\Hilm'$ yields a
  well-defined \Mpn{G_2}equivariant \Mpn{\ast}homomorphism $\rho \colon E \to
  \Adj(\Hilm')$.  Let $F' \defeq \rho(F)$, $\phi' \defeq \rho \circ \phi$.
  Then\/ $(\Hilm', \phi', F')$ is a Kasparov triple for $A,B$.  There is a
  canonical homotopy\/ $(\bar{\Hilm}, \bar{\phi}, \bar{F})$ between\/ $(\Hilm,
  \phi, F)$ and\/ $(\Hilm', \phi', F')$.
  
  If\/ $(\Hilm, \phi, F)$ is an \Mpn{\Hils}special Kasparov triple, then\/
  $(\Hilm', \phi', F')$ and\/ $(\bar{\Hilm}, \bar{\phi}, \bar{F})$ are
  \Mpn{\Hils}special Kasparov triples as well.
\end{lemma}

\begin{proof}
  Since $J \subideal E$ is an ideal, $E(\Hilm') \subset \Hilm'$.  If $T \in
  \Adj(\Hilm)$ satisfies $T(\Hilm') \subset \Hilm'$ and $T^\ast(\Hilm')
  \subset \Hilm'$, then the restriction of~$T$ to~$\Hilm'$ is an adjointable
  operator $\rho(T) \colon \Hilm' \to \Hilm'$.  This yields the desired map
  $\rho \colon E \to \Adj(\Hilm')$.  Since $(\Hilm, \phi, F)$ is a Kasparov
  triple, $J \subset \Comp (\Hilm)$.  We have defined~$\Hilm'$ so that even $J
  \subset \Comp_{\Hilm} (\Hilm')$.  Hence $\rho(J) \subset \Comp(\Hilm')$ by
  Lemma~\ref{lem:embedding_functorial}.  This means that $(\Hilm', \phi', F')$
  is a Kasparov triple.

  Let~$\bar{\Hilm}$ be the Hilbert \Mp{B[0,1],G_2}module $\{ f\in
  \Hilm{}[0,1] \mid f(1) \in \Hilm' \}$.  Define $\bar{F} \in
  \Adj(\bar{\Hilm})$ and $\bar{\phi} \colon A \to \Adj(\bar{\Hilm})$ by
  $(\bar{F}f)(t) \defeq F f(t)$ and $(\bar{\phi}(a) f)(t) \defeq \phi(a)
  f(t)$ for all $a \in A$, $f \in \bar{\Hilm}$, $t \in [0,1]$.  An argument
  similar to the proof that $(\Hilm', \phi', F')$ is a Kasparov triple shows
  that $(\bar{\Hilm}, \bar{\phi}, \bar{F})$ is a Kasparov triple for
  $A,B[0,1]$.  It provides the desired homotopy between $(\Hilm, \phi, F)$
  and $(\Hilm', \phi', F')$.

  Clearly, $(\Hilm', \phi', F')$ and $(\bar{\Hilm}, \bar{\phi}, \bar{F})$ are
  \Mpn{\Hils}special if $(\Hilm, \phi, F)$ is \Mpn{\Hils}special.
\end{proof}

\section{Some universal algebra}
\label{sec:universal}

In this section, we recall the definitions and some elementary properties of
the algebras $\littleQ A$ and $\littleX A$ introduced by
Cuntz~\cite{cuntz:newlook} and Haag~\cite{haag:gradedKK}.  We examine their
relationship to special Kasparov triples and utilize this to describe
$KK^G(A,B)$ as a set of homotopy classes of equivariant homomorphisms.

\subsection{The algebras $\littleX A$, $\midX A$, and $\bigX A$}
\label{sec:X_def}

Let~$A$ be a \Cstar{}algebra.  Define~$\bigX A$ as the universal (unital)
\Cstar{}algebra generated by~$A$ and a symmetry \cite{haag:gradedKK}.  That
is, we have a \Mpn{\ast}homomorphism $j_A \colon A \to \bigX A$ and a
symmetry $F_A \in \bigX A$ such that for all triples $(B,\phi,F)$ consisting
of a unital \Cstar{}algebra~$B$, a \Mpn{\ast}homomorphism $\phi \colon A \to
B$, and a symmetry $F \in B$, there is a unique unital \Mpn{\ast}homomorphism
$(\phi,F)_\ast \colon \bigX A \to B$ satisfying $(\phi,F)_\ast \circ j_A =
\phi$ and $(\phi, F)_\ast (F_A) = F$.

The construction of $\bigX A$ is clearly functorial.  Hence if~$A$ is a
\Mpn{G}\Cstar{}algebra, then there is an induced action of~$G$ on $\bigX A$.
This action is uniquely determined by the requirement that $j_A$ be
\Mpn{G}equivariant and~$F_A$ be \Mpn{G}invariant.  Since non-commutative
polynomials in $j_A(a)$, $a \in A$, and~$F_A$ are dense in $\bigX A$, the
\Mpn{G}action on $\bigX A$ is strongly continuous.  If~$A$ is graded, then we
endow $\bigX A$ with the unique grading~$\tau$ for which~$j_A$ is equivariant
and~$F_A$ is \emph{odd}, that is, $\tau(F_A) = -F_A$.

If $\phi \colon A \to B$ is a \Mpn{G_2}equivariant \Mpn{\ast}homomorphism,
then the induced map $\bigX \phi \colon \bigX A \to \bigX B$ is a
\Mpn{G_2}equivariant \Mpn{\ast}homomorphism as well.

Let $\littleX A \subideal \bigX A$ be the ideal generated by the \emph{graded}
commutators $[j_A(a),F]$ with $a \in A$.  The ideal~$\littleX A$ is
\Mpn{G_2}invariant and essential.  Thus $\bigX A \subset \Mult( \littleX A)$.
The quotient $\bigX A / \littleX A$ is the universal unital \Cstar{}algebra
generated by~$A$ and a symmetry that graded commutes with~$A$.  Thus $\bigX A
/ \littleX A \cong \Cl_1 \hot \Unse{A}$, where~$\Unse{A}$ is the
\Cstar{}algebra obtained by adjoining a unit to~$A$, with $\Unse{A}/ A = \C$.
Let $\midX A \subideal \bigX A$ be the ideal generated by $j_A (A)$.  It
follows that $\midX A / \littleX A \cong \Cl_1 \hot A$, so that we have a
canonical extension of \Mpn{G_2}\Cstar{}algebras
\begin{equation}  \label{eq:ext_xX}
  \littleX A \injto \midX A \prto \Cl_1 \hot A.
\end{equation}
It is shown in the proof of \cite[Theorem 3.6]{haag:gradedKK} that this
extension has a natural---hence \Mpn{G_2}equivariant---completely positive
section.  Roughly speaking, $\midX A$ is the universal \Cstar{}algebra
generated by~$A$ and a symmetry in the \emph{multiplier algebra} $\Mult(\midX
A)$.

Let $A$ and~$B$ be \Mpn{G_2}\Cstar{}algebras.  There is a canonical map $\midX
(A \hot B) \to \midX A \hot B$ that restricts to a map $\littleX (A \hot B)
\to \littleX A \hot B$.  It is defined by the homomorphism $j_A \hot \ID_B
\colon A \hot B \to \midX A \hot B$ and the symmetry $F_A \hot 1 \in
\Mult(\midX A \hot B)$.  For $B = \NBC([0,1])$, we obtain that $\midX$
and~$\littleX$ are \emph{homotopy functors}.  That is, if $f_0,f_1 \colon A
\to A'$ are homotopic, then $\littleX f_0, \littleX f_1 \colon \littleX A \to
\littleX A'$ are homotopic as well.  For $A = \C$, we obtain canonical maps
$\littleX B \to (\littleX \C) \hot B$ and $\midX B \to (\midX \C) \hot B$.
Our next goal is to show that these maps are \Mp{KK}equivalences.  We follow
arguments in the proof of \cite[Proposition~3.8]{haag:gradedKK} in the
non-equivariant case.

\begin{proposition}  \label{pro:midX_homotopy}
  Let~$A$ be a \Mpn{G_2}\Cstar{}algebra.  Then the canonical
  \Mpn{\ast}homomorphism $\ID \hot j_A \colon \Comp(\Z_2\N) A \to
  \Comp(\Z_2\N) \midX A$ is a homotopy equivalence.
\end{proposition}

\begin{proof}
  We call a map of the form $x \mapsto \bigl( \begin{smallmatrix} x & 0 \\ 0 &
  0 \end{smallmatrix} \bigr)$ an \emph{upper left corner embedding}.  We will
  exhibit a canonical homomorphism $f \colon \midX A \to \hat{\Mat}_2 A$ such
  that $f \circ j_A$ and $(\ID_{\hat{\Mat}_2} \hot j_A) \circ f$ are both
  homotopic to the upper left corner embeddings $A \to \hat{\Mat}_2 A$ and
  $\midX A \to \hat{\Mat}_2 \midX A$.  It follows that $\ID_{\Comp(\Z_2\N)}
  \hot f$ is a homotopy inverse for $\ID \hot j_A$.
  
  The homomorphism~$f$ is defined by requiring $f \circ j_A$ to be the upper
  left corner embedding and $f(F_A)$ to be the \emph{standard symmetry}
  \begin{equation}  \label{eq:def_S}
    S \defeq
    \left( \begin{smallmatrix} 0 & 1 \\ 1 & 0 \end{smallmatrix} \right).
  \end{equation}
  By definition, $f \circ j_A \colon A \to \hat{\Mat}_2 A$ is equal to the
  upper left corner embedding.  The symmetries $S$ and $F' \defeq F_A \oplus
  -F_A$ in $\hat{\Mat}_2 \Mult( \midX A)$ anti-commute.  Hence $t \mapsto
  \sqrt{1-t^2} \cdot S + t F'$ is a path of \Mpn{G}equivariant symmetries in
  $\hat{\Mat}_2 \Mult(\midX A)$ connecting them.  This path yields a homotopy
  between $(\ID_{\hat{\Mat}_2} \hot \midX A) \circ f \colon \midX A \to
  \hat{\Mat}_2 \midX A$ and the upper left corner embedding $\midX A \to
  \hat{\Mat}_2 \midX A$.
\end{proof}

Hence the canonical map $\midX A \to (\midX \C) \hot A$ is invertible in
$KK^{G_2} (\midX A, (\midX \C) \hot A)$.

\begin{proposition}  \label{pro:split_littleX}
  Let~$A$ be a separable \Mpn{G_2}\Cstar{}algebra.  Then the canonical map
  $\littleX A \to (\littleX \C) \hot A$ is invertible in $KK^{G_2} (\littleX
  A, (\littleX \C) \hot A)$.
\end{proposition}

\begin{proof}
  This canonical map is part of a morphism of extensions from $\littleX A
  \injto \midX A \prto \Cl_1 \hot A$ to $(\littleX \C) \hot A \to (\midX \C)
  \hot A \prto \Cl_1 \hot A$, where the map $\midX A \to (\midX \C) \hot A$ is
  a \Mp{KK}equivalence by Proposition~\ref{pro:midX_homotopy} and the map
  $\Cl_1 \hot A \to \Cl_1 \hot A$ is the identity map.  Since the two
  extensions have completely positive \Mpn{G_2}equivariant sections, the long
  exact sequences in \Mp{KK}theory are available.  The Five Lemma yields that
  the map $\littleX A \to (\littleX \C) \hot A$ is a \Mp{KK}equivalence as
  well.
\end{proof}

\subsection{The algebras $\littleQ A$ and $\bigQ A$}
\label{sec:Q_def}

Let $\bigQ A \defeq A \ast A$ be the free product of two copies of~$A$
\cite{cuntz:newlook}.  Thus there are two \Mpn{\ast}homomorphisms $\iota^+_A,
\iota^-_A \colon A \to \bigQ A$ such that for any triple $(B,\phi^+,\phi^-)$
consisting of a \Cstar{}algebra~$B$ and a pair of \Mpn{\ast}homomorphisms
$\phi^+,\phi^- \colon A \to B$, there is a unique \Mpn{\ast}homomorphism
$\phi^+ \ast \phi^- \colon \bigQ A \to B$ satisfying $(\phi^+ \ast \phi^-)
\circ \iota^{\pm}_A = \phi^{\pm}$.

Let $\littleQ A \subideal \bigQ A$ be the ideal that is generated by the
differences $\iota^+(a) - \iota^-(a)$ with $a \in A$.  Alternatively, we can
describe~$\littleQ A$ as the kernel of the homomorphism $\ID_{A} \ast \ID_{A}
\colon \bigQ A \to A$.  Thus we obtain an extension of \Cstar{}algebras
$\littleQ A \injto \bigQ A \prto A$.  The \Mpn{\ast}homomorphisms
$\iota^\pm_A \colon A \to \bigQ A$ are sections for $\ID_{A} \ast \ID_{A}$.
There is a natural \Mpn{\ast}homomorphism $\pi_A \defeq (\ID_{A} \ast
0)|_{\littleQ A} \colon \littleQ A \to A$.

If~$A$ is a \Mpn{G_2}\Cstar{}algebra, then there is a unique strongly
continuous \Mpn{G_2}action on~$\bigQ A$ for which the
\Mpn{\ast}homomorphisms~$\iota^{\pm}_A$ are \Mpn{G_2}equivariant.  The
ideal~$\littleQ A$ is \Mpn{G_2}invariant.  The maps $\iota^\pm_A$, $\pi_A$,
and $\ID_{A} \ast \ID_{A}$ above are \Mpn{G_2}equivariant.  The functor $A
\mapsto \littleQ A$ is a homotopy functor.

\begin{proposition}  \label{pro:Q_stable_homotopy}
  Let $A$ and~$B$ be \Mpn{G_2}\Cstar{}algebras.

  Let $\iota_1 \colon A \to A \oplus B$ and $\iota_2 \colon B \to A \oplus B$
  be the standard inclusions.

  The homomorphism $\ID_{\Comp(\N)} \otimes (\iota_1 \ast \iota_2) \colon
  \Comp(\N) (A \ast B) \to \Comp(\N) (A \oplus B)$ is a homotopy equivalence.
  In particular, $\Comp(\N) \bigQ A$ is homotopy equivalent to $\Comp(\N) (A
  \oplus A)$.
\end{proposition}

\begin{proof}
  The stable homotopy inverse for $\iota_1 \ast \iota_2$ is the map $f \colon
  A \oplus B \to \Mat_2(A \ast B)$,
  $$
  f(a,b) \defeq
  \left( \begin{smallmatrix}
    a & 0 \\ 0 & b
  \end{smallmatrix} \right) \qquad \text{for $a \in A$, $b \in B$.}
  $$
  The compositions $f \circ (\iota_1 \ast \iota_2)$ and $\bigl(\ID_{{\Mat_2}}
  \otimes (\iota_1 \ast \iota_2) \bigr) \circ f$ are homotopic to the upper
  left corner embeddings $A \ast B \to \Mat_2(A \ast B)$ and $A \oplus B \to
  \Mat_2(A \oplus B)$ in a natural way.  Roughly speaking, the homotopies
  leave~$a$ fixed and rotate~$b$ to the upper left corner.  Consequently,
  $\ID_{\Comp(\N)} \otimes (\iota_1 \ast \iota_2)$ is a homotopy equivalence.
  The occurring homotopies are natural and therefore \Mpn{G_2}equivariant.
\end{proof}

\subsection{Universal algebras and Kasparov triples}
\label{sec:uni_KK}

\begin{proposition}  \label{pro:KKGs_XQ}
  Let $A$ and~$B$ be \Mpn{\sigma}unital \Mpn{G_2}\Cstar{}algebras and
  let\/~$\Hils$ be a separable \Mpn{G_2}Hilbert space.  There are
  natural bijections
  $$
  KK^G_s (A,B) \cong [\littleX A, \Comp(G_2\N) B],
  \qquad
  KK^G_{s,\Hils} (A,B) \cong [\Comp(\Hils) \littleX A, \Comp(G_2\N) B].
  $$
  If $A$, $B$, and\/~$\Hils$ are trivially graded, then there are natural
  bijections
  $$
  KK^G_s(A,B) \cong [\littleQ A, \Comp(G\N) B],
  \qquad
  KK^G_{s,\Hils} (A,B) \cong [\Comp(\Hils) \littleQ A, \Comp(G\N) B].
  $$
\end{proposition}

All the sets $KK^G_s (A,B)$, $[\littleX A, \Comp(G_2\N) B]$, etc., in the
proposition are functorial for \Mpn{G_2}equivariant \Mpn{\ast}homomorphisms $A
\to A'$, $B' \to B$.  Naturality means that the isomorphisms are compatible
with this functoriality, so that we have isomorphisms of functors, not just of
sets.

\begin{proof}
  Since special Kasparov triples are nothing but \Mpn{\C}special Kasparov
  triples, it suffices to prove the assertions about $KK^G_{s,\Hils}$.  Let
  $T \defeq (\Hilm, \phi, F)$ be an \Mpn{\Hils}special Kasparov triple.  The
  pair $(\phi, F)$ defines a \Mpn{G_2}equivariant \Mpn{\ast}homomorphism
  $(\phi, F)_\ast \colon \bigX A \to \Adj(\Hilm)$ whose restriction to
  $\littleX A$ has values in $\Comp (\Hilm)$.  Hence we get a map $(\phi,
  F)_\ast^{\Hils} \defeq \ID_{\Comp(\Hils)} \hot (\phi, F)_\ast \colon
  \Comp(\Hils) \littleX A \to \Comp(\Hils) \Comp(\Hilm) \cong \Comp(\Hils
  \hot \Hilm)$.  Since the Kasparov triple~$T$ is \Mpn{\Hils}special,
  there is an isometric embedding $\iota \colon \Hils \hot \Hilm \to
  \Hilsg_B$.  Let $\Psi(T) \colon \littleX A \to \Comp(G_2\N) B$ be the
  composition $\Comp(\iota) \circ (\phi, F)_\ast^{\Hils}$.
  
  The homomorphism $\Psi(T)$ is determined uniquely up to homotopy by
  Lemma~\ref{lem:stabilization_homotopy}.  Since we can apply~$\Psi$ to
  \Mpn{\Hils}special homotopies as well, it descends to a map on homotopy
  classes $\Psi \colon KK^G_{s,\Hils} (A, B) \to [\Comp(\Hils) \littleX A,
  \Comp(G_2\N) B]$.  It is straightforward to verify that~$\Psi$ is natural.
  That is, if $f \colon A' \to A$ and $g \colon B \to B'$ are
  \Mpn{G_2}equivariant \Mpn{\ast}homomorphisms, and $T \in KK^G_{s, \Hils}
  (A,B)$, then $\Psi\bigl( f^\ast(T) \bigr) = \Psi(T) \circ
  (\ID_{\Comp(\Hils)} \hot \littleX f)$ and $\Psi\bigl( g_\ast(T) \bigr) =
  (\ID_{\Comp(G_2\N)} \hot g) \circ \Psi(T)$---even if~$g$ is not essential.

  Conversely, let $f \colon \Comp(\Hils) \littleX A \to \Comp(G_2\N)B \cong
  \Comp(\Hilsg_B)$ be a \Mpn{G_2}equivariant \Mpn{\ast}homomorphism.  Let
  $\Hilm_1 = f(\Comp(\Hils) \littleX A) \cdot \Hilsg_B \subset \Hilsg_B$ and
  let $\iota \colon \Hilm_1 \to \Hilsg_B$ be the inclusion mapping.  By
  construction, $f(\Comp(\Hils) \littleX A) \subset \Comp_{\Hilsg_B}
  (\Hilm_1)$.  Hence Lemma~\ref{lem:embedding_functorial} yields $f =
  \Comp(\iota) \circ f_1$ for a \Mpn{G_2}equivariant essential
  \Mpn{\ast}homomorphism $f_1 \colon \Comp(\Hils) \littleX A \to
  \Comp(\Hilm_1)$.
  
  We claim that $\Hilm_1 \cong \Hils \hot \Hilm_2$ and $f_1 \cong
  \ID_{\Comp(\Hils)} \hot f_2$ for a Hilbert \Mp{B,G}module~$\Hilm_2$ and an
  essential \Mpn{G}equivariant \Mpn{\ast}homomorphism $f_2 \colon \littleX A
  \to \Comp(\Hilm_2)$.  This is trivial if $\Hils = \C$.  Consider the dual
  $(\Hils^\ast \hot \littleX A, \psi^\ast)$ of the \Mp{\Comp(\Hils) \littleX
    A, \littleX A, G_2}imprimitivity bimodule $\Hils \hot \littleX A$.  Thus
  $(\Hils \hot \littleX A) \hot_{\psi^\ast} (\Hils^\ast \hot \littleX A) \cong
  \Comp(\Hils) \littleX A$.  Let
  $$
  \Hilm_2 \defeq (\Hils^\ast \hot \littleX A) \hot_{f_1} \Hilm_1
  \qquad \text{and} \qquad
  f_2 \defeq \psi^\ast \hot 1.
  $$
  Since~$\psi^\ast$ is essential, so is~$f_2$.  Since~$f_1$ is essential as
  well, we have $\Hils \hot \Hilm_2 \cong (\Hils \hot \littleX A) \hot_{f_2}
  \Hilm_2 \cong \Hilm_1$.  Under this isomorphism, $f_1$ corresponds to
  $\ID_{\Comp(\Hils)} \hot f_2$.  Since $f_1(\Comp(\Hils) \littleX A) \subset
  \Comp(\Hilm_1)$, it follows that $f_2(\littleX A) \subset \Comp(\Hilm_2)$.
  
  We may extend~$f_2$ to $\bigX A \subset \Mult(\littleX A)$.  By the
  universal property of $\bigX A$, this extension is of the form $(\phi,
  F)_\ast \colon \bigX A \to \Comp(\Hilm_2)$ for some \Mpn{G_2}equivariant
  \Mpn{\ast}homomorphism $\phi \colon A \to \Adj(\Hilm_2)$ and some
  \Mpn{G}invariant symmetry $F \in \Adj(\Hilm_2)$.  The triple $\Psi^{-1}(f)
  \defeq (\Hilm_2, \phi, F)$ is a Kasparov triple because $(\phi,F)_\ast
  (\littleX A) \subset \Comp(\Hilm_2)$.  It is \Mpn{\Hils}special because~$F$
  is a \Mpn{G}equivariant symmetry and $\Hils \hot \Hilm_2 \cong \Hilm_1
  \subset \Hilsg_B$.  Evidently, $\Psi^{-1}$ descends to a map $[\Comp(\Hils)
  \littleX A, \Comp(G_2\N)B] \to KK^G_{s, \Hils} (A,B)$.  By construction,
  $$
  \Comp(\iota) \circ (\ID_{\Comp(\Hils)} \hot (\phi, F)_\ast) =
  \Comp(\iota) \circ (\ID_{\Comp(\Hils)} \hot f_2) =
  \Comp(\iota) \circ f_1 =
  f.
  $$
  That is, $\Psi \circ \Psi^{-1}$ is the identity map on $[\Comp(\Hils)
  \littleX A, \Comp(G_2\N) B]$.
  
  Let $(\Hilm, \phi, F)$ be an \Mpn{\Hils}special Kasparov triple.  Going
  through the above constructions, we find that $\Psi^{-1} \circ \Psi(\Hilm,
  \phi, F)$ is the Kasparov triple that is called $(\Hilm', \phi', F')$ in
  Lemma~\ref{lem:degenerate_homotopy}.  Therefore, $[\Psi^{-1} \circ \Psi
  (\Hilm, \phi, F)] = [(\Hilm, \phi, F)]$ in $KK^G_{s,\Hils} (A,B)$.  The
  proof of the isomorphism $KK^G_{s,\Hils} (A,B) \cong [\Comp(\Hils) \littleX
  A, \Comp(G_2\N) B]$ is finished.
  
  Suppose now that $A$, $B$, and~$\Hils$ are trivially graded.  Let $(\Hilm,
  \phi, F)$ be an \Mpn{\Hils}special Kasparov triple for $A$, $B$.  The even
  and odd part $\Hilm^+$ and~$\Hilm^-$ of~$\Hilm$ are Hilbert \Mp{B,G}modules
  as well.  We may use~$F$ to identify $\Hilm^+ \cong \Hilm^-$.  Then~$F$
  becomes the standard symmetry $S \in \Adj(\Hilm^+ \oplus \Hilm^+)$
  of~\eqref{eq:def_S}.  Since~$A$ is trivially graded, we have $\phi = \phi_+
  \oplus \phi_-$ for certain \Mpn{\ast}homomorphisms $\phi^\pm \colon A \to
  \Adj (\Hilm^+)$.  The condition $[F, \phi(a)] \in \Comp (\Hilm)$ becomes
  $\phi_+(a) - \phi_-(a) \in \Comp (\Hilm^+)$ for all $a \in A$.  Thus
  \Mpn{\Hils}special Kasparov triples correspond bijectively to
  \Mpn{G}equivariant \Mpn{\ast}homomorphisms $f \colon \bigQ A \to
  \Adj(\Hilm)$ with $f(\littleQ A) \subset \Comp(\Hilm)$ and $\Hils \otimes
  \Hilm \subset \Hilsu_A$.  Copying the argument above with $\littleQ A
  \subideal \bigQ A$ instead of $\littleX A \subideal \bigX A$, we obtain the
  desired bijection $KK^G_{s, \Hils} (A, B) \cong [\Comp(\Hils) \littleQ A,
  \Comp(G\N)B]$ if $A$ and~$B$ are trivially graded.
\end{proof}

If $K = \Comp(G\N)$ or $K = \Comp(G_2\N)$, let $[A,B]_K$ be the set of
homotopy classes of \Mpn{G_2}equivariant \Mpn{\ast}homomorphisms from $K \hot
A$ to $K \hot B$.  Let
$$
\littleX_s A \defeq \littleX (\Comp(G_2\N) A)
\qquad \text{and} \qquad
\littleQ_s A \defeq \littleQ (\Comp(G\N) A).
$$

\begin{theorem}  \label{the:KKG_XQ}
  Let~$G$ be a locally compact, \Mpn{\sigma}compact topological group.  Let
  $A$ and~$B$ be \Mpn{\sigma}unital \Mpn{G_2}\Cstar{}algebras.  Let\/
  $\Hils_1$ and\/~$\Hils_2$ be separable \Mpn{G_2}Hilbert spaces.

  There are natural bijections
  \begin{displaymath}
    KK^G(A,B) \cong
    [\Comp(\Hils_1)\, \littleX (\Comp(L^2G \hot \Hils_2)A),\,
    \Comp(G_2\N) B]
    \cong
    [\littleX_s A, B]_{\Comp(G_2\N)}.
  \end{displaymath}
  If $A$, $B$, $\Hils_1$, and\/~$\Hils_2$ are trivially graded, then there are
  natural bijections
  \begin{displaymath}
    KK^G(A,B) \cong
    [\Comp(\Hils_1)\, \littleQ (\Comp(L^2G \otimes \Hils_2) A),\,
    \Comp(G\N) B]
    \cong
    [\littleQ_s A, B]_{\Comp(G\N)}.
  \end{displaymath}
\end{theorem}

The sets $KK^G(A,B)$, etc., occurring in the Theorem are functorial for
equivariant \Mpn{\ast}homomorphisms $A' \to A$, $B \to B'$.  The naturality
of the isomorphisms means that they are compatible with this functoriality.

\begin{proof}
  Since Morita-Rieffel equivalent \Mpn{G_2}\Cstar{}algebras are
  \Mp{KK^G}equivalent, there are natural isomorphisms
  $$
  KK^G(A,B) \cong
  KK^G(\Comp(G) A, B) \cong
  KK^G(\Comp(G) \Comp(\Hils_2) A, B).
  $$
  By Proposition~\ref{pro:stable_aep}, $\Comp(G) A$ has the property AE.
  Hence Proposition~\ref{pro:special_aep} yields $KK^G(\Comp(G) A, B) \cong
  KK^G_s(\Comp(G) A, B) \cong KK^G_{s, \Hils_1} (\Comp(G) A,B)$.  A similar
  statement holds for $\Comp(L^2G \hot \Hils_2) A$ instead of $\Comp(G) A$.
  Therefore, Proposition~\ref{pro:KKGs_XQ} yields the assertions.
\end{proof}

\section{The universal property of equivariant Kasparov theory}
\label{sec:Kasparov_universal}

In this section, we formulate and establish the universal property of
equivariant Kasparov theory for trivially graded separable
\Mpn{G}\Cstar{}algebras.

Let~$\GCalg$ be the category of separable \Mpn{G}\Cstar{}algebras with
\Mpn{G}equivariant \Mpn{\ast}homomorphisms as morphisms.  Let~$\SHo$ be the
\Mp{\Comp(G\N)}stable homotopy category, whose objects are the separable
\Mpn{G}\Cstar{}algebras and whose set of morphisms from~$A$ to~$B$ is $[A,B]_s
\defeq [A,B]_{\Comp(G\N)}$.  Let~$KK^G$ be the category whose objects are the
separable \Mpn{G}\Cstar{}algebras and whose set of morphisms from $A$ to~$B$
is $KK^G(A,B)$.  The Kasparov product yields the composition of morphisms
in~$KK^G$.  We rely on Kasparov's work~\cite{kasparov88:equivariantKK} and
assume that the Kasparov product exists and is associative.  We do not attempt
an alternative definition of the Kasparov product as in~\cite{cuntz:newlook}.
It is clear that~$KK^G$ is an additive category.  There are obvious functors
$\GCalg \to \SHo$ and $\GCalg \to KK^G$.

Let~$\Cat$ be a category.  A functor $F \colon \GCalg \to \Cat$ is called a
\emph{homotopy functor} iff $F(f_0) = F(f_1)$ whenever $f_0$ and~$f_1$ are
\Mpn{G}equivariantly homotopic.

A functor $F \colon \GCalg \to \Cat$ is called \emph{stable} iff the map
$F(\Comp(\Hils) A) \to F(\Comp(\Hils \oplus \Hils') A)$ induced by the
inclusion $\Hils \subset \Hils \oplus \Hils'$ is an isomorphism for all
separable \Mpn{G}Hil\-bert spaces $\Hils,\Hils'$ and all separable
\Mpn{G}\Cstar{}algebras~$A$.

\begin{proposition}  \label{pro:sho_universal}
  The functor $\GCalg \to \SHo$ is a stable homotopy functor.  A functor $F
  \colon \GCalg \to \Cat$ is a stable homotopy functor iff it can be factored
  through the functor $\GCalg \to \SHo$.  This factorization is automatically
  unique.

  In other words, $\SHo$ is the universal stable homotopy functor.
\end{proposition}

\begin{proof}
  It is left to the reader to check that the canonical functor $\GCalg \to
  \SHo$ is a stable homotopy functor.  Thus any functor $\GCalg \to \Cat$ that
  factors through it is a stable homotopy functor as well.

  Conversely, let $F \colon \GCalg \to \Cat$ be a stable homotopy functor.
  Let $\Hils = \C \oplus L^2(G\N)$ and let $j_1^A \colon A \to \Comp(\Hils) A$
  and $j_2^A \colon \Comp(G\N) A \to \Comp(\Hils) A$ be the canonical
  inclusions.  Since~$F$ is stable, $F(j_1^A)$ and $F(j_2^A)$ are
  isomorphisms.  Thus $\sigma_A \defeq F(j_2^A)^{-1} \circ F(j_1^A)$ is a
  natural isomorphism $F(A) \congto F(\Comp(G\N) A)$.  Define
  $$
  F_\ast \colon [A,B]_s \to \Mor_{\Cat}\bigl( F(A), F(B) \bigr),
  \qquad
  F_\ast [\phi] \defeq \sigma_B^{-1} \circ F(\phi) \circ \sigma_A.
  $$
  It is left to the reader to check that this defines a functor $F_\ast \colon
  \SHo \to \Cat$ that extends~$F$ and that the functor~$F_\ast$ is determined
  uniquely.
\end{proof}

\begin{remark}
  A homotopy functor $F \colon \GCalg \to \Cat$ is stable iff $F(A) \cong
  F(\Comp(G\N) A)$ naturally.  The proof of
  Proposition~\ref{pro:sho_universal} shows that a natural isomorphism $F(A)
  \cong F(\Comp(G\N) A)$ allows us to factor~$F$ through~$\SHo$.  Our
  definition of a stable homotopy functor is equivalent to the definitions
  in~\cite{guentner-higson-trout} and in~\cite{thomsen:KKGuniversal}.
\end{remark}

A functor $F \colon \GCalg \to \Cat$ into an additive category~$\Cat$ is
called \emph{split exact} iff $\bigl( F(i), F(s) \bigr) \colon F(A) \oplus
F(C) \to F(B)$ is an isomorphism for all extensions
$$
0 \to
A \xrightarrow{i}
B \xrightarrow{p}
C \to
0
$$
of \Mpn{G}\Cstar{}algebras that split by a \Mpn{G}equivariant
\Mpn{\ast}homomorphism $s \colon C \to B$.

\begin{proposition}  \label{pro:KKG_properties}
  The canonical functor $\GCalg \to KK^G$ is a split exact stable homotopy
  functor.
\end{proposition}

\begin{proof}
  Clearly, $KK^G$ is a stable homotopy functor.  Split exactness is a
  straightforward consequence of the associativity of the Kasparov product.
  The argument in \cite[Proposition~2.1]{cuntz:newlook} carries over without
  change.
\end{proof}

Since $A \mapsto \littleQ A$ is a homotopy functor, $A \mapsto \littleQ_s A$
descends to a functor from $\SHo$ to itself.  The map $\pi_{\Comp(G\N) A}
\colon \littleQ_s A \to \Comp(G\N) A$ gives rise to a natural morphism
$\pi^s_A \in [\littleQ_s A, A]_s$.

\begin{lemma}  \label{lem:split_q}
  Let $F \colon \GCalg \to \Cat$ be a split exact stable homotopy functor and
  let $F_\ast \colon \SHo \to \Cat$ be the unique extension of~$F$.  Then
  $F_\ast (\pi^s_A)$ is invertible for all~$A$.
\end{lemma}

\begin{proof}
  Split exactness applied to the extension $A \injto A \oplus B \prto B$
  yields that the canonical map $F(A \oplus B) \to F(A) \oplus F(B)$ is an
  isomorphism.  That is, $F$ is \emph{additive}.
  Proposition~\ref{pro:Q_stable_homotopy} yields $F(A \ast B) \cong F(A)
  \oplus F(B)$.  Split exactness applied to the extension $\littleQ A \injto
  \bigQ A \prto A$ implies that $F(\pi_A) \colon F(\littleQ A) \to F(A)$ is
  an isomorphism for all~$A$.  This implies that $F_\ast(\pi^s_A)$ is an
  isomorphism as well.
\end{proof}

By Proposition~\ref{pro:KKG_properties} and
Proposition~\ref{pro:sho_universal}, the canonical functor $\GCalg \to KK^G$
factors through a functor $\natural \colon \SHo \to KK^G$.
Lemma~\ref{lem:split_q} implies that $\natural(\pi^s_A) \in KK^G(\littleQ_s
A, A)$ is invertible for all~$A$.

\begin{theorem}  \label{the:KKG_doubleQ}
  Let $A$ and~$B$ be separable \Mpn{G}\Cstar{}algebras.  The map
  $$
  [\littleQ_s A, \littleQ_s B]_s \to
  KK^G(A,B),
  \qquad
  f \mapsto
  \natural(\pi^s_B) \circ \natural(f) \circ \natural(\pi^s_A)^{-1},
  $$
  is a natural isomorphism.  Hence the Kasparov product on $KK^G$ corresponds
  to the composition of homomorphisms.
\end{theorem}

\begin{proof}
  Since $\pi^s_B$ induces an isomorphism $KK^G(A, \littleQ_s B) \cong KK^G(A,
  B)$, it suffices to verify that the isomorphism $[\littleQ_s A, B]_s \to
  KK^G(A,B)$ of Theorem~\ref{the:KKG_XQ} is given by $f \mapsto \natural(f)
  \circ \natural(\pi^s_A)^{-1}$.  By naturality, it suffices to check this for
  the identity map in $[\littleQ_s A, \littleQ_s A]_s$.  Composing with the
  invertible element~$\pi^s_A$, we can reduce the theorem to the following
  claim: The isomorphism of Theorem~\ref{the:KKG_XQ} maps $\pi^s_A \in
  [\littleQ_s A, A]_s$ to the unit in $KK^G (A,A)$, represented by the
  Kasparov triple $(A, \ID_{A}, 0)$.  The proof of this claim is made somewhat
  messy by stabilizations, but otherwise straightforward.  Therefore, we omit
  it.
\end{proof}

\begin{theorem}  \label{the:universal_KKG}
  The functor $\GCalg \to KK^G$ is the universal split exact stable homotopy
  functor in the following sense.  An additive functor $F \colon \GCalg \to
  \Cat$ into an additive category~$\Cat$ can be extended to a functor $F_\ast
  \colon KK^G \to \Cat$ iff it is a split exact stable homotopy functor.  The
  extension is necessarily unique.
\end{theorem}

\begin{proof}
  Let $F \colon \GCalg \to \Cat$ be a split exact stable homotopy functor.  By
  Proposition~\ref{pro:sho_universal}, we may assume that~$F$ is a functor $F
  \colon \SHo \to \Cat$.  Split exactness implies that $F(\pi^s_A)$ is an
  isomorphism for all~$A$.  If $f \in [\littleQ_s A, \littleQ_s B]_s$, define
  $F_\ast(f) \defeq F(\pi_s^B) \circ F(f) \circ F(\pi_s^A)^{-1}$.  By
  Theorem~\ref{the:KKG_doubleQ}, this yields a functor $KK^G \to \Cat$.
  Evidently, this is the unique functor extending~$F$.  It is clear that any
  additive functor that factors through $KK^G$ is a split exact stable
  homotopy functor.
\end{proof}

\section{The case of graded algebras}
\label{sec:graded}

Following Haag~\cite{haag:gradedKK}, we write $\Ex^G (A,B) \defeq KK^{G_2}
(A,B)$ for the \Mpn{G_2}equivariant \Mp{KK}theory for trivially graded
algebras.  We show $KK^G (A,B) \cong \Ex^G (\hat{S} \hot A, B)$ and describe
the Kasparov product in $KK^G$ in terms of the product in $\Ex^G$.

We redefine $KK^G$ to be the category whose objects are the \Mpn{\Ztwo}graded
separable \Mpn{G}\Cstar{}algebras and whose set of morphisms from $A$ to~$B$
is $KK^G (A,B)$.  Let $\GCalg[G_2]$ be the category of separable
\Mpn{G_2}\Cstar{}algebras and let $\SHo[G_2]$ be the \Mp{\Comp(G_2\N)}stable
homotopy category, as defined in the previous section.  We redefine
$\littleQ_s A \defeq \littleQ (\Comp(G_2\N) A)$, so that $\Ex^G (A, B) \cong
[\littleQ_s A, B]_s$ by Theorem~\ref{the:KKG_XQ}.

The canonical functor $\GCalg[G_2] \to KK^G$ is still a split exact stable
homotopy functor.  By Theorem~\ref{the:universal_KKG}, we may extend it to a
functor $\alpha \colon \Ex^G (A,B) \to KK^G (A,B)$.  The functor~$\alpha$ can
be computed as follows.  As in \cite[p.~15]{haag:gradedKK}, the
\Mpn{\ast}homomorphism $\iota^+ \oplus \iota^- \colon A \to \hat{\Mat}_2
(\bigQ A)$ and the symmetry~$S$ of~\eqref{eq:def_S} yield a canonical map
$$
\alpha_0 \defeq (\iota^+ \oplus \iota^-, S)_\ast \colon
\littleX A \to \hat{\Mat}_2 ( \littleQ A)
$$
We view~$\alpha_0$ as an element of $[\littleX A, \littleQ A]_s$.
Replacing~$A$ by $\Comp(G_2\N) A$, we obtain $\alpha_0 \in [\littleX_s A,
\littleQ_s A]_s \cong KK^G (A, \littleQ_s A)$ by Theorem~\ref{the:KKG_XQ}.

\begin{lemma}  \label{lem:alpha0}
  $\alpha_0 \in KK^G (A, \littleQ_s A)$ is the inverse of\/ $\pi^s_A \in
  [\littleQ_s A, A]_s$.
\end{lemma}

\begin{proof}
  Lemma~\ref{lem:split_q} implies that the image of~$\pi^s_A$ in $KK^G
  (\littleQ_s A, A)$ is invertible.  It remains to prove that $(\pi^s_A)_\ast
  (\alpha_0)$ is the identity element of $KK^G (A,A)$.
  
  We may suppose $\Comp(G_2\N) A \cong A$, so that we may omit the
  stabilizations and work with the map $\alpha_0 \colon \littleX A \to
  \hat{\Mat}_2 (\littleQ A)$.  It corresponds to the Kasparov triple
  $(\littleQ A \oplus (\littleQ A)^\opp, \iota^+ \oplus \iota^-, S)$.  Since
  $\pi_A \circ \iota^+ = \ID_A$, $\pi_A \circ \iota^- = 0$, we have
  $$
  (\pi_A)_\ast (\alpha_0) =
  (A \oplus A^\opp, \ID_A \oplus 0, S).
  $$
  The right hand side represents the identity element of $KK^G (A,A)$.
\end{proof}

\begin{corollary}  \label{cor:alpha}
  Let $A$ and~$B$ be separable \Mpn{G_2}\Cstar{}algebras.  Using the
  isomorphisms of Theorem~\ref{the:KKG_XQ}, we obtain a map
  $$
  [\littleQ_s A, B]_s \cong
  \Ex^G (A,B) \xrightarrow{\alpha}
  KK^G (A,B)  \cong
  [\littleX_s A, B]_s.
  $$
  This map is equal to composition with $[\alpha_0] \in [\littleX_s A,
  \littleQ_s A]_s$.
\end{corollary}

\begin{proof}
  Let $f \in [\littleQ_s A,B]_s$, then the image of~$f$ in $KK^G(A,B)$ is
  $f_\ast (\pi^s_A)^{-1} = f_\ast [\alpha_0]$.  This is mapped to $f \circ
  [\alpha_0] \in [\littleX_s A, B]_s$.
\end{proof}

There is a canonical Kasparov triple $(\littleX A, j_A, F_A)$ for $A, \littleX
A$.  Replacing~$A$ by $\Comp(G_2\N) A$, we obtain a canonical element $i_A \in
KK^G (A, \littleX_s A)$.  The isomorphism $KK^G (A, \littleX_s A) \to
[\littleX_s A, \littleX_s A]_s$ maps~$i_A$ to the identity map.  The
naturality of the isomorphism $[\littleX_s A, B] \to KK^G(A,B)$ of
Theorem~\ref{the:KKG_XQ} implies that it maps $f \mapsto f_\ast(i_A)$ for all
$f \in [\littleX_s A,B]_s$.

\begin{lemma}  \label{lem:ExlittleX_KK}
  Let $A$ and~$B$ be separable \Mpn{G_2}\Cstar{}algebras.  The canonical map
  $$
  \Ex^G (\littleX_s A, B) \xrightarrow{\alpha}
  KK^G(\littleX_s A, B)   \xrightarrow{i_A^\ast}
  KK^G(A,B)
  $$
  is an isomorphism.
\end{lemma}

\begin{proof}
  Theorem~\ref{the:KKG_XQ} yields canonical isomorphisms $\Ex^G (\littleX_s A,
  B) \cong [\littleQ_s \littleX_s A, B]_s$ and $KK^G (A, B) \cong [\littleX_s
  A, B]_s$.  We are going to show that $\pi \defeq \pi^s_{\littleX_s A} \colon
  \littleQ_s \littleX_s A \to \littleX_s A$ is invertible in $\SHo[G_2]$.
  Therefore, $[\littleQ_s \littleX_s A, B]_s \cong [\littleX_s A, B]_s$.  It
  is straightforward to show that the corresponding isomorphism $\Ex^G
  (\littleX_s A, B) \cong KK^G (A,B)$ is equal to the map in the statement
  Lemma~\ref{lem:ExlittleX_KK}.
  
  The homotopy inverse for~$\pi$ is constructed as a Kasparov product.  Let $i
  = i_A \in KK^G(A, \littleX_s A)$ be as above.  Let $j \in KK^G(\littleX_s A,
  \littleQ_s \littleX_s A)$ be the inverse of~$\pi$.  Let $h \in KK^G(A,
  \littleQ_s \littleX_s A) \cong [\littleX_s A, \littleQ_s \littleX_s A]_s$ be
  the Kasparov product of $i$ and~$j$.  The associativity of the Kasparov
  product implies $\pi \circ h = i$ in $KK^G(A, \littleX_s A) = [\littleX_s A,
  \littleX_s A]_s$.  Since~$\pi$ is invertible in $\Ex^G$, composition
  with~$\pi$ is an isomorphism
  $$
  \Ex^G (\littleX_s A, \littleQ_s \littleX_s A) \congto
  \Ex^G (\littleX_s A, \littleX_s A).
  $$
  Hence the equality $\pi \circ h \circ \pi = \pi$ in $[\littleQ_s
  \littleX_s A, \littleX_s A]_s$ implies $h \circ \pi = \ID$ in $[\littleQ_s
  \littleX_s A, \littleQ_s \littleX_s A]_s$.  Thus~$h$ is inverse to~$\pi$ in
  $\SHo[G_2]$.
\end{proof}

Let~$\hat{S}$ be the algebra $\CVI(\R)$ graded by $\tau f(x) = f(-x)$ for all
$x \in \R$, $f \in \CVI(\R)$ and with trivial \Mpn{G}action.  It is shown in
the proof of \cite[Proposition~3.8]{haag:gradedKK} that $\littleX \C \cong
\hat{\Mat}_2 \hat{S}$, so that $\hat{S}$ and $\littleX \C$ are Morita-Rieffel
equivalent.  Together with Proposition~\ref{pro:split_littleX}, we obtain a
canonical isomorphism in $\Ex^G (\littleX_s A, \hat{S} \hot A)$.

Let $e \in KK^G(\C, \hat{S}) \cong [\littleX_s \C, \hat{S}]_s$ be represented
by the isomorphism $\littleX \C \to \hat{\Mat}_2 \hat{S}$.  It is easy to
verify that~$e$ is homotopic to the Kasparov triple $(\hat{S}, 1,
x/\sqrt{1+x^2})$, where $1 \colon \C \to \Adj(\hat{S}) \cong \BC(\R)$ is the
unique unital \Mpn{\ast}homomorphism and $x/ \sqrt{1+x^2}$ denotes the bounded
function $x \mapsto x/ \sqrt{1+x^2}$ on~$\R$.

\begin{theorem}  \label{the:graded_equi_KK}
  Let~$G$ be a locally compact \Mpn{\sigma}compact topological group and let
  $A$ and~$B$ be separable \Mpn{G_2}\Cstar{}algebras.
  The composition
  $$
  \sigma \colon 
  \Ex^G (\hat{S} \hot A, B) \xrightarrow{\alpha}
  KK^G (\hat{S} \hot A, B) \xrightarrow{(e \hot \ID_A)^\ast}
  KK^G (A,B)
  $$
  is an isomorphism.  Here $(e \hot \ID_A)^\ast$ denotes the Kasparov
  product with the exterior product $e \hot \ID_A \in KK^G (A, \hat{S} \hot
  A)$.
\end{theorem}

\begin{proof}
  The isomorphism $KK^G(A, \littleX_s A) \cong KK^G(A, (\littleX \C) \hot A)$
  induced by the canonical map $\littleX_s A \to (\littleX \C) \hot
  \Comp(G_2\N) A$ maps~$i_A$ to the exterior product $i_\C \hot \ID_A$.  Hence
  the isomorphism $KK^G (A, \littleX_s A) \to KK^G(A, \hat{S} \hot A)$
  maps~$i_A$ to $e \hot \ID_A$.  If we compose the isomorphism $\Ex^G
  (\littleX_s A, B) \to KK^G (A,B)$ of Lemma~\ref{lem:ExlittleX_KK} with the
  the isomorphism $\Ex^G(\hat{S} \hot A, B) \to \Ex^G (\littleX_s A,B)$
  induced by the \Mp{\Ex^G}equivalence $\hat{S} \hot A \to \littleX_s A$, we
  obtain that~$\sigma$ is an isomorphism.
\end{proof}

We have to compute the exterior product $e \hot e \in KK^G(\C, \hat{S} \hot
\hat{S})$.  Since the \Mpn{G}action on~$\littleX \C$ and $\hat{S}$ is trivial,
we may forget about the \Mpn{G}actions.  Therefore, we briefly resort to the
case of trivial~$G$.  Theorem~\ref{the:KKG_XQ} implies
$$
KK (\C, B) = [\littleX \C, \Comp(\Z_2\N) B] \cong [\hat{S}, \Comp(\Z_2\N) B].
$$
We claim that $e \hot e \in KK (\C, \hat{S} \hot \hat{S})$ belongs to
the homomorphism $\hat{S} \to \hat{S} \hot \hat{S}$ that is called~$l$
by Haag \cite[p.~87]{haag:algebraic} and~$\Delta$ by Higson and
Kasparov~\cite{higson-kasparov:actHilbert}.

To verify this elementary claim, it is convenient to describe $KK (A,B)$ by
unbounded operators following Baaj and Julg~\cite{baaj-julg:unbounded} because
in this picture exterior products are straightforward to compute.  The
unbounded picture of $KK(\C, B)$ is also nicely related to the isomorphism
$KK(\C, B) \cong [\hat{S}, \Comp(\Z_2\N) B]$.  The essential, grading
preserving \Mpn{\ast}homomorphisms $\hat{S} \to \Adj(\Hilm)$ correspond
bijectively to odd, self-adjoint, possibly unbounded multipliers of~$\Hilm$
via $f \mapsto f(\ID_\R)$ for $f \colon \hat{S} \to \Adj(\Hilm)$.  Since~$e$
belongs to the unbounded multiplier~$\ID_\R$ of~$\hat{S}$, the exterior
product $e \hot e$ belongs to the unbounded multiplier $\ID_\R \hot 1 + 1 \hot
\ID_\R$ of $\hat{S} \hot \hat{S}$.  Thus $e \hot e$ is represented by the
map~$\Delta$ of~\cite{higson-kasparov:actHilbert}.  It is easy to check that
the concrete formula for~$l$ in~\cite{haag:gradedKK} yields nothing
but~$\Delta$.

\begin{theorem}  \label{the:graded_product}
  Let $A$, $B$, and~$C$ be \Mpn{G_2}\Cstar{}algebras and let $x \in \Ex^G
  (\hat{S} \hot A,B)$, $y \in \Ex^G (\hat{S} \hot B,C)$.  The Kasparov product
  of $\sigma(y) \in KK^G (B,C)$ and $\sigma(x) \in KK^G (A,B)$ is mapped
  by~$\sigma^{-1}$ to the composition
  $$
  \hat{S} \hot A \xrightarrow{\Delta \hot \ID_A}
  \hat{S} \hot \hat{S} \hot A \xrightarrow{\ID_{\hat{S}} \hot y}
  \hat{S} \hot B \xrightarrow{x}
  C
  $$
  in $\Ex^G$.
\end{theorem}

\begin{proof}
  Recall the definition of~$\sigma$ in Theorem~\ref{the:graded_equi_KK} and
  that~$\alpha$ is multiplicative.  Moreover, it is easy to check
  that~$\alpha$ is compatible with exterior products, so that
  $\alpha(\ID_{\hat{S}} \hot x) \cong \ID_{\hat{S}} \hot \alpha(x)$.  Hence we
  compute
  \begin{multline*}
    \sigma(y) \circ \sigma(x)
    =
    \alpha(y) \circ (e \hot \ID_B) \circ \alpha(x) \circ (e \hot \ID_A)
    \\ =
    \alpha(y) \circ \bigl( e \hot \alpha(x) \bigr) \circ (e \hot \ID_A)
    =
    \alpha(y) \circ \bigl(\ID_{\hat{S}} \hot \alpha(x) \bigr) \circ
    (e \hot \ID_{\hat{S} \hot A}) \circ (e \hot \ID_A)
    \\ =
    \alpha(y) \circ \alpha(\ID_{\hat{S}} \hot x) \circ (e \hot e \hot \ID_A)
    =
    \sigma\bigl(y \circ (\ID_{\hat{S}} \hot x) \circ
    (\Delta \hot \ID_A) \bigr).
  \end{multline*}
  We used that the Kasparov product is compatible with exterior products.
\end{proof}

\section{Proper actions and square-integrable Hilbert modules}
\label{sec:proper}

Exel~\cite{exel:MoritaSpectral} and Rieffel~\cite{rieffel:pre1} define
the concept of a proper action of a locally compact group on a
\Cstar{}algebra.  Furthermore, Rieffel relates proper \Mpn{G}actions on the
algebra $\Comp(\Hils)$ to square-integrable representations of~$G$.  It is
very illuminating to consider also square-integrable group actions on Hilbert
modules.  The main result is that a countably generated Hilbert
\Mp{A,G}module is square-integrable iff it is a direct summand of~$\Hilsu_A$.
We conclude that proper algebras have property AE.

Concerning questions of properness, we may ignore gradings whenever
convenient.  Since the group~$\Ztwo$ is compact, a \Mpn{G_2}\Cstar{}algebra
is proper iff it is proper as a \Mpn{G}\Cstar{}algebra.

Let~$A$ be a \Mpn{G}\Cstar{}algebra and let~$\Hilm$ be a Hilbert
\Mp{A,G}module.  We denote the \Mpn{G}actions on $A$ and~$\Hilm$ by $\alpha$
and~$\gamma$, respectively.  We frequently view~$A$ as a right Hilbert
\Mp{A,G}module.  Let $(K_n)_{n \in \N}$ be a sequence of compact subsets
of~$G$ such that $K_{n+1}$ is a neighborhood of~$K_n$ for all~$n$ and $G =
\bigcup K_n$.  Let $(\kappa_n)_{n \in \N}$ be an increasing sequence of
functions $\kappa_n \colon G \to [0,1]$ with $\kappa_n|_{K_n} = 1$ and
$\kappa_n|_{G \setminus K_{n+1}} = 0$.

A continuous function $f \colon G \to A$ is called \emph{square-integrable}
iff the sequence $\int_G f^\ast(g) f(g) \kappa_n(g) \,dg$ is a norm Cauchy
sequence in~$A$.  Equivalently, the sequence of integrals $\int_{K_n}
f^\ast(g) f(g) \,dg$ is norm convergent.  Observe that these sequences are
increasing sequences of positive elements and that the notion of
square-integrability does not depend on the choice of the sets~$K_n$ or the
functions~$\kappa_n$.

It is easy to check that~$f$ is square-integrable iff the sequence $(f \cdot
\kappa_n)_{n \in \N}$ is a Cauchy sequence with respect to the norm $\norm{h}
\defeq \norm{\int_G h^\ast(g) h(g) \,dg}^{1/2}$ on $\CC(G,A)$.  Since the
completion of $\CC(G,A)$ with respect to this norm is precisely $L^2(G,A)$, we
can view square-integrable continuous functions as elements of $L^2(G,A)$.

If $\xi, \eta \in \Hilm$, then we define the \emph{coefficient function}
$c_{\xi\eta} \colon G \to A$ by
$$
c_{\xi\eta}(g) \defeq \5{\gamma_g(\xi)} {\eta}_A
\qquad \text{for all $g \in G$.}
$$
In the special case $\Hilm = A$, we have $c_{ab}(g) \defeq \alpha_g(a)^\ast
b$.

We call $\xi \in \Hilm$ \emph{square-integrable} iff the
function~$c_{\xi\eta}$ is square-integrable for all $\eta \in \Hilm$.  The
Hilbert module~$\Hilm$ is called \emph{square-integrable} iff the set of
square-integrable elements is dense in~$\Hilm$.  A \Mp{G}\Cstar{}algebra~$A$
is called \emph{proper} iff it is square-integrable as a right Hilbert
\Mp{A,G}module.  Let $A_+ \subset A$ be the cone of positive elements.  We
call $a \in A_+$ \emph{integrable} iff~$a^{1/2}$ is square-integrable.

By definition, $a \in A_+$ is integrable iff the integrals $\int_{K_n} b^\ast
\alpha_g(a) b \,dg$ form a Cauchy sequence with respect to the
norm topology for all $b \in A$.  Moreover, $a \in A$ is square-integrable
iff $a a^\ast$ is integrable.  Hence~$A$ is proper iff the set of integrable
elements is dense in~$A_+$.  The above definition of properness is equivalent
to Rieffel's definition in~\cite{rieffel:pre1} and thus also to Exel's
definition in~\cite{exel:MoritaSpectral}.

\begin{lemma}  \label{lem:square_integrable}
  Let~$\Hilm$ be a Hilbert \Mp{A,G}module and let $\xi, \eta, \zeta \in
  \Hilm$.
  \begin{enumerate}[(i)]
  \item If~$\xi$ is square-integrable, then the map
    $$
    \Gamma_\xi \colon \Hilm \to L^2(G,A), \qquad
    \Gamma_\xi(\eta) \defeq c_{\xi\eta},
    $$
    is adjointable.  The adjoint\/ $\Gamma_\xi^\ast \colon L^2 (G,A) \to
    \Hilm$ satisfies
    \begin{equation}  \label{eq:gamma_ast_def}
      \Gamma_\xi^\ast(f) \defeq \int_G \gamma_g(\xi) \cdot f(g) \,dg
      \qquad \text{for all $f \in \CC(G, A)$.}
    \end{equation}

  \item The operators\/ $\Gamma_\xi$ and\/~$\Gamma_\xi^\ast$ are
    \Mpn{G}equivariant.
    
  \item The closure of the range of\/~$\Gamma_\xi^\ast$ is the smallest
    \Mpn{G}invariant Hilbert submodule of~$\Hilm$ containing~$\xi$.  In
    particular, $\xi$ is contained in the closure of\/ $\Ran \Gamma_\xi^\ast$.

  \item If $\xi$ and~$\zeta$ are square-integrable, then the sequence
    $$
    \int_G \gamma_g(\xi) \cdot \5{\gamma_g(\zeta)} {\eta}_A \,\kappa_n(g)
    \,dg, \qquad n \in \N,
    $$
    in~$\Hilm$ is norm convergent.  Its limit is\/ $\Gamma_\xi^\ast
    \Gamma_\zeta (\eta)$.
 
  \item $\xi$ is square-integrable iff\/ $\ket{\xi} \bra{\xi} \in
    \Comp(\Hilm)$ is integrable.

  \end{enumerate}
\end{lemma}

\begin{proof}
  The Banach-Steinhaus theorem yields that~$\Gamma_\xi$ is bounded.  We can
  define an operator $\Gamma_\xi^\ast \colon \CC(G,A) \to \Hilm$
  by~\eqref{eq:gamma_ast_def}.  For $f \in \CC(G,A)$, we compute
  \begin{equation}  \label{eq:Txi_adjoint}
    \5{\Gamma_\xi^\ast(f)}{\eta}_A =
    \int_G \5{\gamma_g(\xi) f(g)}{\eta}_A \,dg =
    \int_G f(g)^\ast \cdot c_{\xi\eta}(g) \,dg =
    \5{f}{\Gamma_\xi (\eta)}_A.
  \end{equation}
  Hence $\norm{ \5{\Gamma_\xi^\ast(f)}{\eta}_A } \le \norm{f}_{L^2(G,A)}
  \norm{\eta} \cdot \norm{\Gamma_\xi}$.  Since~$\eta$ is arbitrary, it follows
  that $\norm{\Gamma_\xi^\ast(f)} \le \norm{f}_{L^2(G,A)} \cdot
  \norm{\Gamma_\xi}$.  Thus we may extend~$\Gamma_\xi^\ast$ to $L^2(G,A)$.
  Equation~\eqref{eq:Txi_adjoint} shows that~$\Gamma_\xi^\ast$ is adjoint
  to~$\Gamma_\xi$.

  Straightforward computations show that~$\Gamma_\xi$ and $\Gamma_\xi^\ast$
  are \Mpn{G}equivariant.
  
  Assertion~(iii) follows easily once we know that~$\xi$ is contained in the
  closed range of~$\Gamma_\xi^\ast$.  Choose $\epsilon >0$.  There is $u \in
  A$ with $0 \le u \le 1$ and $\norm{ \xi\cdot u - \xi} \le \epsilon/2$.
  There is a compact neighborhood~$U$ of $1 \in G$ with $\norm{\gamma_g(\xi) -
  \xi} < \epsilon/2$ for all $g \in U$.  Let $f \colon G \to \R_+$ be a
  continuous function with support~$U$ and $\int_G f(g) \,dg = 1$.  Then
  $\norm{\Gamma_\xi^\ast (f \otimes u) - \xi} \le \epsilon$.  Hence~$\xi$ is
  contained in the closure of $\Ran \Gamma_\xi^\ast$.

  We compute
  \begin{equation}  \label{eq:Gxi_Gxi_ad}
    \Gamma_\xi^\ast \Gamma_\zeta(\eta) =
    \Gamma_\xi^\ast (c_{\zeta\eta}) =
    \lim_{n\to \infty} \Gamma_\xi^\ast (c_{\zeta\eta}\kappa_n) =
    \lim_{n\to \infty} \int_G \gamma_g(\xi) \5{\gamma_g(\zeta)}{\eta}_A
    \, \kappa_n(g) \,dg.
  \end{equation}
  The boundedness of~$\Gamma_\xi^\ast$ implies that the sequence is norm
  convergent.

  Equation~\eqref{eq:Gxi_Gxi_ad} implies that the sequence
  $$
  I_n \defeq \int_G \ket{\gamma_g(\xi)} \bra{\gamma_g(\xi)}
  \, \kappa_n(g) \,dg
  = \int_G \gamma_g(\ket{\xi} \bra{\xi}) \, \kappa_n(g) \,dg
  \in \Comp(\Hilm)
  $$
  is bounded and converges strongly (that is, pointwise on~$\Hilm$) towards
  $\Gamma_\xi^\ast \Gamma_\xi$.  Therefore, the sequences $(I_n \cdot T)$ and
  $(T \cdot I_n)$ converge in norm for all $T \in \Comp(\Hilm)$.  This means
  that $\ket{\xi} \bra{\xi} \in \Comp(\Hilm)$ is integrable.  Conversely, if
  $\ket{\xi} \bra{\xi} \in \Comp(\Hilm)$ is integrable, then the sequence
  $\5{\eta} {I_n(\eta)}_A$ is norm convergent for all $\eta \in \Comp(\Hilm)
  \cdot \Hilm = \Hilm$.  Since
  $$
  \5{\eta} {I_n(\eta)}_A =
  \int_G \5{\gamma_g(\xi)}{\eta}_A^\ast \5{\gamma_g(\xi)}{\eta}_A
  \,\kappa_n(g) \,dg =
  \int_G c_{\xi\eta}(g)^\ast c_{\xi\eta}(g) \, \kappa_n(g) \,dg,
  $$
  this means that~$\xi$ is square-integrable.
\end{proof}

\begin{remark}
  If $\Gamma_\xi^\ast \colon \CC(G,A) \to \Hilm$ extends to an
  \emph{adjointable} map $L^2(G,A) \to \Hilm$, then~$\xi$ is
  square-integrable.
  
  The map~$\Gamma_\xi^\ast$ extends to a bounded operator $L^2(G,A) \to \Hilm$
  if and only if the function $g \mapsto c_{\xi\eta}(g)^\ast c_{\xi\eta}(g)$
  is \emph{order-integrable} in the sense of Rieffel~\cite{rieffel:pre1} for
  all $\eta \in \Hilm$.  Hence it may happen that~$\Gamma_\xi^\ast$ extends to
  a bounded operator on $L^2(G,A)$ that is not adjointable.
\end{remark}

\begin{proposition}  \label{pro:si_proper}
  Let~$A$ be a \Mpn{G}\Cstar{}algebra and let $\Hilm$ be a Hilbert
  \Mp{A,G}module.  Then~$\Hilm$ is square-integrable iff\/ $\Comp(\Hilm)$ is
  proper.
\end{proposition}

\begin{proof}
  If~$\Hilm$ is square-integrable, then the linear span of the integrable
  elements is dense in $\Comp(\Hilm)$ by
  Lemma~\ref{lem:square_integrable}.(v).  Therefore, $\Comp(\Hilm)$ is
  proper.  Conversely, assume that $\Comp(\Hilm)$ is proper.  Let $T \in
  \Comp(\Hilm)_+$ be square-integrable.  If $\xi \in \Hilm$, then $\ket{T\xi}
  \bra{T\xi} = T \ket{\xi} \bra{\xi} T^\ast \le \norm{\xi}^2 TT^\ast$ is
  integrable because $TT^\ast$ is integrable and the set of integrable
  elements is a hereditary cone in~$\Comp(\Hilm)_+$~\cite{rieffel:pre1}.
  Hence $T\xi \in \Hilm$ is square-integrable by
  Lemma~\ref{lem:square_integrable}.(v).  Since $\Comp(\Hilm)$ is proper, the
  set of elements of~$\Hilm$ of the form $T\xi$ with square-integrable $T \in
  \Comp(\Hilm)$ is dense in~$\Hilm$.  Thus~$\Hilm$ is square-integrable.
\end{proof}

\begin{proposition}  \label{pro:proper_functorial}
  Let $A$ and~$B$ be \Mpn{G}\Cstar{}algebras, let~$\Hilm$ be a Hilbert
  \Mp{B,G}module, and let $\phi \colon A \to \Adj(\Hilm)$ be an essential
  \Mpn{G}equivariant \Mpn{\ast}homomorphism.

  If~$A$ is proper, then~$\Hilm$ is square-integrable.
\end{proposition}

\begin{proof}
  Identify $\Adj(\Hilm) \cong \Mult\bigl( \Comp(\Hilm) \bigr)$.  By
  \cite[Theorem~5.3]{rieffel:pre1}, we conclude that $\Comp(\Hilm)$ is
  proper.  Thus~$\Hilm$ is square-integrable by
  Proposition~\ref{pro:si_proper}.
\end{proof}

\begin{theorem}  \label{the:square_integrable}
  Let~$A$ be a \Mpn{G}\Cstar{}algebra and let~$\Hilm$ be a countably
  generated Hilbert \Mp{A,G}module.  Then the following assertions are
  equivalent:
  \begin{enumerate}[(i)]
  \item $\Hilm$ is square-integrable;
  \item $\Comp(\Hilm)$ is proper;
  \item there is a \Mpn{G}equivariant unitary $\Hilm \oplus \Hilsu_A \cong
    \Hilsu_A$;
  \item $\Hilm$ is a direct summand of\/~$\Hilsu_A$.
  \end{enumerate}
\end{theorem}

\begin{proof}
  Proposition~\ref{pro:si_proper} asserts that (i) and~(ii) are equivalent.
  It is trivial that (iii) implies~(iv).  It remains to show that~(iv)
  implies~(i) and that~(i) implies~(iii).
  
  We prove that~(iv) implies~(i).  It is straightforward to show that
  $\Comp(L^2G)$ is a proper \Mpn{G}\Cstar{}algebra.  Equivalently, $L^2G$ is a
  square-integrable \Mpn{G}Hilbert space.  Let~$\Hilm[F]$ be an arbitrary
  Hilbert \Mp{A,G}module.  The canonical \Mpn{\ast}homomorphism $\Comp(L^2G)
  \to \Adj(L^2G \otimes \Hilm[F])$, $T \mapsto T \otimes 1$, is essential and
  \Mpn{G}equivariant.  Hence $L^2(G, \Hilm[F])$ is square-integrable by
  Proposition~\ref{pro:proper_functorial}.  Especially, $\Hilsu_A$ is
  square-integrable.  A direct summand of a square-integrable Hilbert module
  is square-integrable as well because the projection onto the direct summand
  maps square-integrable elements to square-integrable elements.  Hence any
  direct summand of~$\Hilsu_A$ is square-integrable.  That is, (iv)
  implies~(i).
  
  The proof that~(i) implies~(iii) is very similar to the proof of the
  stabilization theorem by Mingo and
  Phillips~\cite{mingo-phillips:triviality}.  Suppose that~$\Hilm$ is
  square-integrable.  Hence there is a sequence $(\xi_n)_{n \in \N}$ of
  square-integrable elements of~$\Hilm$, whose linear span is dense
  in~$\Hilm$.  Let $\Gamma_n \defeq \Gamma_{\xi_n}$ be as in
  Lemma~\ref{lem:square_integrable}.  We may assume that $\norm{\Gamma_n}
  \le 1$ for all $n \in \N$ and that each~$\xi_n$ is repeated infinitely
  often.  An element of~$\Hilsu_A$ can be viewed as a sequence $(f_n)_{n \in
    \N}$ with $f_n \in L^2(G,A)$.  We formally write $\sum f_n \delta_n$ for
  this sequence.  Define an adjointable operator $T \colon \Hilsu_A \to
  \Hilm \oplus \Hilsu_A$ by
  \begin{gather*}
  T\biggl( \sum_{n=1}^\infty f_n \delta_n \biggr) \defeq
  \sum_{n=1}^\infty 2^{-n} \Gamma_n^\ast(f_n) \oplus
  \sum_{n=1}^\infty 4^{-n} f_n \delta_n,
  \\
  T^\ast|_{\Hilm} (\eta) \defeq
  \sum_{n=1}^\infty 2^{-n} \Gamma_n(\eta) \delta_n,
  \qquad
  T^\ast|_{\Hilsu_A} \biggl( \sum_{n=1}^\infty f_n \delta_n \biggr) \defeq
  \biggl( \sum_{n=1}^\infty 4^{-n} f_n \delta_n \biggr).
  \end{gather*}
  Lemma~\ref{lem:square_integrable}.(ii) implies that $T$ is
  \Mpn{G}equivariant.  Evidently, $T^\ast$ has dense range.

  We claim that~$T$ has dense range as well.  Let $\Hilm[F] \subset \Hilm
  \oplus \Hilsu_A$ be the closure of the range of~$T$.  Let $f \in L^2(G,A)$.
  Since each~$\Gamma_n^\ast$ is repeated infinitely often, we have
  $\Gamma_n^\ast(f) \oplus 2^{-k}f \delta_k \in \Ran T$ for infinitely many
  $k \in \N$.  Hence $\Gamma_n^\ast(f) \oplus 0 \in \Hilm[F]$ for all $f \in
  L^2(G,A)$.  By Lemma~\ref{lem:square_integrable}.(iii), this implies $\xi_n
  \oplus 0 \in \Hilm[F]$ for all~$n$ and hence $\Hilm \subset \Hilm[F]$.
  Finally, we get $0 \oplus f \delta_n \in \Hilm[F]$ for all $f \in L^2(G,A)$
  and thus $\Hilm[F] = \Hilm \oplus \Hilsu_A$.
  
  Since both $T$ and~$T^\ast$ have dense range, the composition $T^\ast T$ has
  dense range.  Thus $\abs{T} \defeq (T^\ast T)^{1/2}$ has dense range because
  $\abs{T} (\Hilm) \supset \abs{T} (\abs{T} \Hilm) = T^\ast T (\Hilm)$.  Since
  $\5{\abs{T}\eta}{\abs{T}\eta}_A = \5{T^\ast T\eta}{\eta}_A =
  \5{T\eta}{T\eta}_A$, the formula $U(\abs{T} \eta) \defeq T\eta$ well-defines
  an isometry~$U$ from $\Ran \abs{T}$ onto $\Ran T$.  Extending~$U$
  continuously, we obtain the desired unitary $U \colon \Hilsu_A \to \Hilm
  \oplus \Hilsu_A$.  Since~$T$ is \Mpn{G}equivariant, so is~$U$.
\end{proof}

Thomsen~\cite{thomsen:extensions} calls a \Mpn{G}\Cstar{}algebra~$A$
\emph{\Mpn{K}proper} iff $\Hilm \oplus \Hilsu_A \cong \Hilsu_A$ for all
Hilbert \Mp{A,G}bimodules~$\Hilm$.  Theorem~\ref{the:square_integrable}
implies that~$A$ is \Mpn{K}proper in Thomsen's sense iff all
\Mpn{G}\Cstar{}algebras that are Morita-Rieffel equivalent to~$A$ are proper
in our sense.  For instance, the algebra $\Comp(G)$ is not \Mpn{K}proper,
unless~$G$ is compact, because it is Morita-Rieffel equivalent to the improper
\Mpn{G}\Cstar{}algebra~$\C$.

\begin{proposition}  \label{pro:proper_aep}
  All \Mp{\sigma}unital proper \Mpn{G_2}\Cstar{}algebras have property AE.
\end{proposition}

\begin{proof}
  Let $A$ and~$B$ be \Mp{\sigma}unital \Mpn{G_2}\Cstar{}algebras and let
  $(\Hilm, \phi, F)$ be an essential Kasparov triple for $A,B$.  Suppose
  that~$A$ is proper.  Proposition~\ref{pro:proper_functorial} implies
  that~$\Hilm$ is square-integrable.  Hence $\Hilm \subset \Hilsg_B$ by
  Theorem~\ref{the:square_integrable}.
  
  Let~$\Hilsg_A^\ast$ be the imprimitivity bimodule dual to~$\Hilsg_A$.  As
  remarked in Section~\ref{sec:imprimitivity}, we have $\Hilsg_A^\ast =
  \Comp(\Hilsg_A, A)$ with a canonical Hilbert \Mp{A, \Comp(\Hilsg_A),
  G_2}bimodule structure.  Since~$A$ is proper and \Mpn{\sigma}unital,
  Theorem~\ref{the:square_integrable} yields a \Mpn{G_2}equivariant,
  adjointable isometry $T \colon A \to \Hilsg_A$.  Composition with~$T$ gives
  rise to an adjointable isometry $T_\ast \colon \Comp(\Hilsg_A, A) \to
  \Comp(\Hilsg_A, \Hilsg_A)$.  Thus $\Hilsg_A^\ast$ is a direct summand in
  $\Comp(\Hilsg_A)$.

  Lemma~\ref{lem:equi_connection} yields a \Mpn{G}equivariant
  \Mpn{F}connection~$\bar{F}$ on $L^2(G_2, \Hilm)^\infty = \Hilsg_A \hot_\phi
  \Hilm$.  If we view~$\bar{F}$ as an operator on $\Comp(\Hilsg_A)
  \hot_{\Comp(\Hilsg_A)} L^2(G_2, \Hilm)^\infty \cong L^2(G_2, \Hilm)^\infty
  $, then we obtain an \Mpn{\bar{F}}connection.  Hence the compression
  $$
  F' \defeq
  (T_\ast \hot_{\Comp(\Hilsg_A)} 1)^\ast \cdot
  \bar{F} \cdot
  (T_\ast \hot_{\Comp(\Hilsg_A)} 1)
  $$
  of~$\bar{F}$ to $\Hilsg_A^\ast \hot_{\Comp(\Hilsg_A)} \Hilsg_A \hot_\phi
  \Hilm \cong A \hot_\phi \Hilm \cong \Hilm$ is an \Mpn{\bar{F}}connection as
  well.  By \cite[18.3.4.f]{blackadar98:Ktheory}, $F'$ is an
  \Mpn{F}connection.  Another \Mpn{F}connection is~$F$ itself.  Therefore,
  $F'$ is a compact perturbation of~$F$.  Since~$\bar{F}$ and~$T$ are
  \Mpn{G}equivariant, so is~$F'$.
\end{proof}

\begin{theorem}  \label{the:proper_KK}
  Let $A$ and~$B$ be \Mpn{\sigma}unital \Mpn{G_2}\Cstar{}algebras.  If~$A$ is
  proper, then
  $$
  KK^G (A,B) \cong
  KK^G_s (A,B) \cong
  [\littleX A, \Comp(G_2\N) B] \cong
  [\littleX A, B]_{\Comp(G_2\N)}.
  $$
  If~$A$ is proper and $A$ and~$B$ are trivially graded, then
  $$
  KK^G (A,B) \cong
  KK^G_s (A,B) \cong
  [\littleQ A, \Comp(G\N) B] \cong
  [\littleQ A, B]_{\Comp(G\N)}.
  $$
\end{theorem}

\begin{proof}
  By Proposition~\ref{pro:proper_aep}, $A$ has property AE.  Hence the
  assertions follow from Proposition~\ref{pro:special_aep} and
  Proposition~\ref{pro:KKGs_XQ}.
\end{proof}

\section{Other equivariant theories}
\label{sec:generalize}

The arguments above generalize to other more general versions of equivariant
Kasparov theory.  First of all, we did not even care to specify whether we
work with complex or real \Cstar{}algebras: The theory above goes through in
both cases.  In the real case, we only have to interpret~$\C$ as the
algebra~$\R$ of real numbers everywhere above.  Hence $KK^G$ is the universal
split exact stable homotopy functor also for real \Cstar{}algebras.  We may
also treat ``real'' \Cstar{}algebras as in~\cite{kasparov:KK}.

Our results carry over to Kasparov's functor $\RKK^G$ with some obvious
changes.  We first define the functors $\RKK^G$ and $RKK^G$.

Let~$G$ be, as usual, a locally compact \Mpn{\sigma}compact topological group
and let~$X$ be a locally compact \Mpn{\sigma}compact \Mpn{G}space.  An
\emph{\Mp{X\rtimes G_2}\Cstar{}algebra} is a \Mpn{G_2}\Cstar{}algebra~$A$
equipped with a \Mpn{G_2}equivariant \emph{essential} \Mpn{\ast}homomorphism
from $\CVI(X)$ into the \emph{center} of $\Mult(A)$.  Let $A,B$ be
\Mp{X\rtimes G_2}\Cstar{}algebras.  A \Mpn{\ast}homomorphism $\phi \colon A
\to \Mult(B)$ is called \emph{\Mp{X \rtimes G_2}equivariant} iff it is
\Mpn{G_2}equivariant and in addition satisfies $\phi(f \cdot a) = f \cdot
\phi(a)$ for all $f \in \CVI(X)$, $a \in A$.

If~$B$ is an \Mpn{X \rtimes G_2}\Cstar{}algebra and~$\Hilm$ is a Hilbert
\Mpn{B,G}module, then $\Comp(\Hilm)$ is an \Mpn{X\rtimes G_2}\Cstar{}algebra
as well.  The homomorphism $\CVI(X) \to \Mult\bigl( \Comp(\Hilm) \bigr) =
\Adj(\Hilm)$ is defined by $f \cdot (\xi \cdot b) \defeq \xi \cdot (f \cdot
b)$ for all $f \in \CVI(X)$, $\xi \in \Hilm$, $b \in B$.  By the Cohen-Hewitt
factorization theorem, all elements of~$\Hilm$ are of the form $\xi \cdot b$
for suitable $\xi \in \Hilm$, $b \in B$.  Computing the inner products $\5{f
  \cdot (\xi \cdot b)}{\eta \cdot c}$ shows that $f \cdot (\xi \cdot b)$ is
well-defined and defines a \Mpn{\ast}homomorphism from $\CVI(X)$ into the
center of $\Adj(\Hilm)$.

Kasparov~\cite{kasparov88:equivariantKK} defines the functor $\RKK^G(X; A,B)$
for \Mp{X\rtimes G_2}\Cstar{}algebras using Kasparov triples $(\Hilm, \phi,
F)$ for $A,B$ with the additional assumption that $\phi \colon A \to
\Adj(\Hilm)$ be \Mp{X\rtimes G_2}equivariant.  If $A$ and~$B$ are just
\Mpn{G_2}\Cstar{}algebras, then he puts $RKK^G(X; A,B) \defeq \RKK^G\bigl(X;
\CVI(X,A), \CVI(X,B) \bigr)$.  Hence $RKK^G$ is a special case of $\RKK^G$.

The presence of the central homomorphism $\CVI(X) \to \Adj(\Hilm)$ creates no
problems in Sections \ref{sec:equi_connections} and~\ref{sec:embed_Hilbert}.
In Section~\ref{sec:universal}, we have to modify the definitions of the
universal algebras $\bigX A$ and $\bigQ A$ because, as they are defined there,
they are not \Mp{X\rtimes G_2}\Cstar{}algebras.  Thus we replace them by
algebras with analogous universal properties in the category of \Mp{X\rtimes
G_2}\Cstar{}algebras.  This amounts to dividing out the relations
$\iota^\pm(fa) \iota^\pm(b) = \iota^\pm(a) \iota^\pm(fb)$ for all $a,b \in A$,
$f \in \CVI(X)$ in $\littleQ A$ and the relation $j_A(fa) F_A j_A(b) = j_A(a)
F_A j_A(fb)$ for all $a,b \in A$, $f \in \CVI(X)$ in $\littleX A$.  The
quotients by the ideals generated by these relations carry a canonical
\Mp{X\rtimes G_2}\Cstar{}algebra structure.  Once this modification is made,
the results of Sections \ref{sec:universal}, \ref{sec:Kasparov_universal},
and~\ref{sec:graded} carry over without change.  Especially, $\RKK^G(X; A,B)$
is the universal split exact stable homotopy functor for trivially graded
separable \Mp{X\rtimes G}\Cstar{}algebras.  Theorem~\ref{the:proper_KK}
remains valid as well.

The functor $\RKK^G(X; A,B)$ is a special case of the equivariant Kasparov
theory for groupoids developed in~\cite{legall:KK_groupoids}.  I expect that
the arguments above carry over to the case of locally compact groupoids with
Haar system.  However, I have not checked the details.  Some work has to be
done to carry over the proof of Lemma~\ref{lem:equi_connection}.  Moreover, to
carry over the results of Section~\ref{sec:proper}, one first has to define
properness in the sense of Rieffel for actions of groupoids.

\bibliographystyle{plain}
\bibliography{kk}

\end{document}